\documentclass[aap,MSNbibl,dvips]{arximspdf}

\doi{10.1214/13-AAP971} 
\volume{24}
\issue{5}
\pubyear{2014}
\firstpage{2033}
\lastpage{2069}

\makeatletter

\def\cal{\mathcal}
\newcommand{\eqref}[1]{(\ref{#1})}

\newtheorem{theorem}{Theorem}[section]

\newproclaim{definition}{Definition}[section]
\newproclaim{example}{Example}[section]
\newtheorem{proposition}{Proposition}[section]
\newproclaim{condition}{Condition}[section]
\newproclaim{assumption}{Assumption}[section]

\newproclaim{remark}{Remark}[section]

\newcommand{\Bf}{\mathfrak{B}}
\newcommand{\E}{\mathbb{E}}
\newcommand{\R}{\mathbb{R}}
\newcommand{\p}{\mathbb{P}}
\newcommand{\N}{\mathbb{N}}
\newcommand{\I}{\mathcal{I}}
\newcommand{\Ir}{\mathbb{I}}
\newcommand{\sig}{\sigma}
\newcommand{\eps}{\varepsilon}

\def\mathnormal{}
\newcommand{\al}{\alpha}
\newcommand{\s}{\sigma}
\newcommand{\Gam}{\mathnormal{\Gamma}}
\newcommand{\Del}{\mathnormal{\Delta}}
\newcommand{\Th}{\mathnormal{\Theta}}
\newcommand{\La}{\mathnormal{\Lambda}}
\newcommand{\Ph}{\mathnormal{\Phi}}
\newcommand{\Ps}{\mathnormal{\Psi}}
\newcommand{\Ups}{\mathnormal{\Upsilon}}
\newcommand{\Om}{\mathnormal{\Omega}}

\newcommand{\EE}{\mathbb{E}}
\newcommand{\PP}{\mathbb{P}}
\newcommand{\bph}{\bolds{\varphi}}
\newcommand{\bxi}{\bolds{\xi}}
\newcommand{\bzeta}{\bolds{\zeta}}
\newcommand{\bw}{\mathbf{w}}
\newcommand{\brho}{\bolds{\rho}}

\newcommand{\calA}{{\cal A}}
\newcommand{\calC}{{\cal C}}
\newcommand{\calD}{{\cal D}}
\newcommand{\calE}{{\cal E}}
\newcommand{\calF}{{\cal F}}
\newcommand{\calI}{{\cal I}}
\newcommand{\calP}{{\cal P}}

\newcommand{\oo}{\bar}
\newcommand{\w}{\wedge}
\newcommand{\iy}{\infty}
\newcommand{\osc}{\operatorname{osc}}

\newcommand{\A}{{\cal A}}
\newcommand{\B}{{\cal B}}
\newcommand{\IA}{\mathit{IA}}
\newcommand{\ST}{\mathit{ST}}
\makeatother

\begin{document}
\begin{frontmatter}

\title{Control of the multiclass $G/G/1$ queue in the moderate deviation
regime\thanksref{T1}}
\thankstext{T1}{Supported in part by the ISF (Grant 1349/08), the
US-Israel BSF (Grant 2008466),
and the Technion fund for promotion of research.}
\runtitle{Control in the moderate deviation regime}

\begin{aug}
\author[A]{\fnms{Rami} \snm{Atar}\corref{}\ead[label=e1]{atar@ee.technion.ac.il}}
\and
\author[A]{\fnms{Anup} \snm{Biswas}}
\runauthor{R. Atar and A. Biswas}
\affiliation{Technion--Israel Institute of Technology}
\address[A]{Department of Electrical Engineering\\
Technion--Israel Institute of Technology\\
Haifa 32000\\
Israel\\
\printead{e1}} 
\end{aug}

\received{\smonth{4} \syear{2012}}
\revised{\smonth{9} \syear{2013}}

%
\begin{abstract}
A multi-class single-server system with general service time
distributions is studied in a moderate deviation heavy traffic regime.
In the scaling limit, an optimal control problem associated with the model
is shown to be governed by a differential game that can be explicitly solved.
While the characterization of the limit by a differential game
is akin to results at the large deviation scale,
the analysis of the problem is closely related to the much studied area
of control in heavy traffic at the diffusion scale.
\end{abstract}

%

\end{frontmatter}

\section{Introduction}\label{sec1}

Models of controlled queueing systems have been studied under various scaling
limits. These include heavy traffic diffusion approximations, which are
based on
the central limit theorem (see \cite{chen-yao,bell-will} and
references therein)
and large deviation (LD) asymptotics; see, for example,
\cite{ata-dup-shw,AGS1} and references therein. To the best of our knowledge,
the intermediate, moderate deviation (MD) scale has not been considered
before in relation to
controlled queueing systems.
In this paper we consider the multi-class $G/G/1$ model in a heavy
traffic MD regime
with a risk-sensitive type cost of a general form,
characterize its asymptotic behavior in terms of a differential game
(DG), and solve the game.
In a special but important case, we also identify a simple policy that
is asymptotically optimal (AO).
The treatment in the MD regime
shares important characteristics with both asymptotic regimes alluded
to above.
It is similar to analogous results in the LD regime, in that the limit
behavior is indeed
governed by a DG. The DG itself is closely related to Brownian control
problems (BCP)
that arise in diffusion approximations. In particular,
the solution method by which BCP are transformed into problems involving
the so-called workload process, turns out to be useful for solving
these DG as well.

Treatments of queueing models in the MD regime without dynamic control
aspects include the following.
In \cite{Puhal-Whitt}, Puhalskii and Whitt prove LD and MD principles
for renewal processes.
Puhalskii \cite{Puhal-1999} establishes LD and MD principles for
queue length and waiting time processes for the single server queue and
for single class
queueing networks in heavy traffic (Puhalskii refers to this regime as
\textit{near} heavy traffic,
to emphasize that the deviations from critical load are at a larger
scale than under standard
heavy traffic; we will use the term heavy traffic in this paper).
Majewski \cite{Majewski06} treats feedforward multi-class network
models with
priority.
Wischik \cite{Wischik} (see also \cite{Ganesh})
illuminates on various links between results on queueing problems in LD
and MD regimes, as well as
similarities between MD and diffusion scale results, particularly the
validity of
results such as the snapshot principle and state space collapse.
Based on these similarities he
conjectures that the well-established dynamic control theory for
heavy traffic diffusion approximations
should have a parallel at the MD scale. Our treatment certainly
confirms this expectation, at least
for the model under investigation.
Cruise \cite{cruise} considers LD and MD as a part of a broader
parametrization framework
for studying queueing systems.

In the model under consideration (see the next section for a complete
description),
customers of $I$ different classes arrive at the system following
renewal processes
and are enqueued in buffers, one for each class. A server, that may
offer simultaneous
service to the various classes, divides its effort among the (at most) $I$
customers waiting at the
head of the line of each buffer. The service time distributions depend
on the class.
The problem is to control these fractions of effort so as to minimize a cost.
MD scaling is obtained by considering a sequence $b_n$, where
$b_n\to\iy,\sqrt{n}/b_n\to\infty$. The arrival and service time
scales are set
proportional to a large parameter $n$, with possible correction of
order $b_n\sqrt n$.
Denoting by $X^n_i(t)$, the number of class-$i$ jobs in the $n$th
system at time $t$,
a scaled version is given by $\tilde X^n=(b_n\sqrt{n})^{-1}X^n$.
Moreover, a critical load condition is assumed, namely that
the limiting traffic intensity is one. The cost is given by
\[
\frac{1}{b_n^2}\log\E\bigl\{e^{b_n^2[\int_0^Th(\tilde X^n(t))\,dt+g(\tilde
X^n(T))]}\bigr\},
\]
where $T>0$, and $h$ and $g$ are given functions.

This type of cost is called \textit{risk-sensitive}; see the book by
Whittle \cite{whittle-book}.
The optimal control formulation of a dynamical system with small noise
goes back to Fleming \cite{fle-1}, who studies the associated
Hamilton--Jacobi equations.
The connection of risk-sensitive cost to DG was made by Jacobson \cite{Jac-1}.
The study of risk-sensitive control via LD theory and the formulation
of the corresponding
maximum principle are due to Whittle \cite{whittle-1}.
Various aspects of this approach have been\vadjust{\goodbreak} studied for controlled
stochastic differential equations,
for example, \cite{dup-mce,fle-mce,fle-sou}.
For queueing networks, risk sensitive control in the LD regime has been studied
in \cite{Dup-1,ata-dup-shw,AGS1}.
Operating a queueing system so as to avoid large queue length or waiting
time is important in
practice, for preventing buffer overflow and assuring quality of service.
A risk-sensitive criterion penalizes such events heavily, and thus
provides a natural way to
address these considerations.
Further motivation for this formulation is that
the solution automatically leads to robustness properties of the
policy; see Dupuis et al. \cite{DJP}.
Note that working in MD scale leads to some additional desired
robustness properties.
Namely, since the rate function in this case typically depends only
on first and second moments of the underlying
primitives, the characteristics of the problem are insensitive to
distributional perturbations which preserve these moments.
The price paid for working in MD scale is that a critical load
condition has to be assumed
for the problem to be meaningful (as it is in diffusion approximations
but not in LD analysis).

The DG governing the limit behavior can be solved explicitly, a fact that
not only is useful in characterizing the limit in a concrete way, but
also turns out to be of crucial importance when proving the convergence.
To describe the game (see Section~\ref{sec2} for the precise definition),
consider the dynamics
\[
\varphi(t)=x+yt+\int_0^t\bigl(\tilde
\lambda(s)-\tilde\mu(s)\bigr)\,ds+\eta (t)\in\R_+^I.
\]
Here $x$ is an initial condition, $y$ is a term capturing the order
$b_n\sqrt n$ time scale correction
alluded to above and $\tilde\lambda$ and $\tilde\mu$ represent
perturbations at scale
$b_n/\sqrt n$ of arrival and service rates, respectively.
These are functions mapping $[0,T]\to\R_+^I$, controlled by player 1.
Next, $\eta\dvtx [0,\iy)\to\R_+^I$ is a function whose formal derivative
represents deviations at scale $b_n/\sqrt n$
of the fraction of effort dedicated by the server to each class.
This function is controlled by player 2 and is considered admissible if:
(a) for all $t$, $\varphi(t)\in\R_+^I$, (b)
$\theta\cdot\eta(0)\ge0$ and (c) $\theta\cdot\eta$ is
nondecreasing, where
$\theta=(\frac{1}{\mu_1},\ldots,\frac{1}{\mu_I})$ is what is
often called the \textit{workload vector}
in the heavy traffic literature.
The cost, which player 1 (resp., 2) attempts to maximize (minimize) is
given by
%
\begin{equation}
\label{62} \int_0^Th\bigl(\varphi(s)
\bigr)\,ds+g\bigl(\varphi(T)\bigr)-\int_0^T\sum
\bigl[a_i\tilde \lambda_i(s)^2+b_i
\tilde\mu_i(s)^2\bigr]\,ds,
\end{equation}
where $a_i$ and $b_i$ are positive constants.

It is instructive to compare this to the game obtained under LD scaling.
The form presented here corresponds to the multiclass $M/M/1$ model, following
\cite{AGS1} (the setting there includes multiple, heterogenous servers,
but the presentation here is specialized to the case of a single server).
One considers
\[
\varphi=\Gam[\psi],\qquad \psi(t)=x+\int_0^t\bigl(
\bar\lambda (s)-u(s)\bullet\bar\mu(s)\bigr)\,ds,
\]
where $\Gam$ is the Skorohod map with normal reflection on the
boundary of
the positive orthant, $\bar\lambda$ and $\bar\mu$ are functions\vadjust{\goodbreak}
$[0,T]\to[0,\iy)^I$,
representing perturbations at the LD scale, and controlled by a
maximizing player;
$u\dvtx [0,T]\to S$ where $S=\{s\in[0,1]^I\dvtx \sum s_i=1\}$ is controlled by
minimizing player
representing fraction of effort per class, and
$\bullet$ denotes the entrywise product of two vectors of the same dimension.
The cost here takes the form
%
\begin{equation}
\label{61} \int_0^Th\bigl(\varphi(s)
\bigr)\,ds+g\bigl(\varphi(T)\bigr) -\int_0^T
\bigl[1\cdot l\bigl(\bar\lambda(s)\bigr)+u(s)\cdot\hat l\bigl(\bar\mu(s)\bigr)
\bigr]\,ds,
\end{equation}
where $l$ and $\hat l$ represent LD cost associated with atypical
behavior; see \cite{AGS1} for more details.
The paper \cite{AGS1} provides a characterization of the game's value
in terms of a Hamilton--Jacobi--Isaacs (HJI) equation. However,
it is not known if the game can be solved explicitly.
In contrast, the game associated with MD turns out to be explicitly
solvable, as we show in this
paper. The reason for this is that
while in the LD game the last term of the cost \eqref{61} involves
both $(\bar\lambda,\bar\mu)$ and $u$,
the corresponding term in \eqref{62}
involves only $(\tilde\lambda,\tilde\mu)$, not $\eta$. Hence this
term plays no role when one computes the optimal response $\eta$
to a given $(\tilde\lambda,\tilde\mu)$ [it does when one optimizes
over $(\tilde\lambda,\tilde\mu)$].
This optimal response is computed via projecting the dynamics in the
direction of the workload vector,
and using minimality considerations of
the one-dimensional Skorohod problem. In fact, the optimal response
$\eta$ to $(\tilde\lambda,\tilde\mu)$
is precisely the one that arises in the diffusion scale analysis of the model,
used there to map the Brownian motion term to the optimal control for
the BCP.
Thus the link to diffusion approximations is strong.

In \cite{AGS1} (following the technique of Atar, Dupuis and Shwartz
\cite{ata-dup-shw}),
the convergence is proved by establishing upper and lower bounds on the
limiting risk-sensitive control
problem's value in terms of the lower and, respectively, upper values
of the DG.
The existence of a limit is then argued via uniqueness of solutions to
the HJI equation satisfied by both values. The arrival and service
are assumed to follow Poisson processes and the convergence proof uses
the form of the Markovian generator and martingale inequalities related
to it.
Since in the MD regime the performance depends only on the first two
moments of the primitives,
these moments carry all relevant information regarding the limit (under
tail assumptions), and so in
this paper we aim at general arrival and service processes.
As a result, the tools based on the Markovian formulation mentioned
above cannot be used.
The approach we take uses completely different considerations.
The asymptotic behavior of the risk-sensitive control problem is
estimated, above and below,
directly by the DG lower value (the corresponding upper value is not
dealt with at all in this paper).
This is made possible thanks to the explicit solvability of the game.
More precisely, the arguments by which the game's optimal strategy is
found,\vadjust{\goodbreak}
including the workload formulation and
the minimality property associated with the Skorohod map,
give rise, when applied to the control problem, to the lower bound.
The proof of the upper bound is by construction
of a particular control which again uses the solution of the game and
its properties.
Note that this approach eliminates the need for any PDE
analysis.\looseness=-1

The control that is constructed in the proof of the upper bound is too
complicated for practical implementation. However, in the case where
$h$ and $g$ are linear
(see Section~\ref{sec5} for the precise linearity condition), a simple
solution to the DG is available,
in the form of a fixed priority policy according to the well-known
$c\mu$ rule.\vadjust{\goodbreak} As our final result shows,
applying a priority policy in the queueing model, according to the same
order of customer classes,
is AO in this case.

To summarize the main contribution of the paper, we have:
\begin{itemize}
\item
provided
the first treatment of a queueing control problem at the MD scale,
\item identified and solved the DG governing the scaling limit for
quite a general setting and
\item proved AO of a simple policy in the linear case.
\end{itemize}
The following conclusions stem from this work:
\begin{itemize} 
\item
Techniques such as the equivalent workload
formulation, which have proven powerful for control problems at the
diffusion scale,
are useful at the MD scale. They are likely to be applicable
in far greater generality than the present setting.
\item
Although control problems at
MD and LD scales are both motivated by similar rationale,
MD is evidently more tractable for the model under consideration, and
potentially
this is true in greater generality.
\end{itemize}

We will use the following notation.
For a positive integer $k$ and $a,b\in\R^k$, $a\cdot b$ denotes the usual
scalar product, while $\|\cdot\|$ denotes Euclidean norm. We denote
$[0, \iy)$ by $\R_+$.
For $T>0$ and a function $f\dvtx [0,T]\to\R^k$,
let $\|f\|^*_t=\sup_{s\in[0,t]}\|f(s)\|$, $t\in[0,T]$. When $k=1$,
we write $|f|^*_t$ for $\|f\|^*_t$.
We sometimes write $\|f\|^*$ for $\|f\|^*_T$ when there is no ambiguity
about $T$.
Denote by $\calC([0,T],\R^k)$ and $\calD([0,T], \R^k)$ the
spaces of continuous functions $[0,T]\to\R^k$ and, respectively,
functions that are
right-continuous with finite left limits (RCLL). Endow the space $\calD
([0,T],\R^k)$ with the
$J_1$ metric,
defined as\looseness=-1
%
\begin{eqnarray}
\label{74} d\bigl(\varphi,\varphi^\prime\bigr)=\inf
_{f\in\Ups} \Bigl(\|f\|^\circ\vee \sup_{[0, T]}
\bigl\|\varphi(t)-\varphi^\prime\bigl(f(t)\bigr)\bigr\| \Bigr),
\nonumber
\\[-8pt]
\\[-8pt]
\eqntext{\varphi,
\varphi^\prime\in\calD\bigl([0, T], \R^k\bigr),}
\end{eqnarray}\looseness=0
where $\Ups$ is the set of strictly increasing, continuous functions
from $[0, T]$ onto itself,
and
%
\begin{equation}
\label{76} \| f\|^\circ=\sup_{0\leq s<t\leq T} \biggl|\log
\frac
{f(t)-f(s)}{t-s} \biggr|.
\end{equation}
As is well known \cite{Bill}, $\calD([0, T], \R^k)$ is a Polish
space under this metric.

The organization of the paper is as follows. The next section
introduces the model and an associated
differential game and states the main result. In Section~\ref{sec3} we find a
solution to the game
and describe properties of it that are useful in the sequel.
Section~\ref{sec4} gives the proof of the main theorem.
In Section~\ref{sec5} we discuss the case of linear cost and identify an AO policy.
Finally, the \hyperref[app]{Appendix} gives the proof of a proposition stated in
Section~\ref{sec2}.\vadjust{\goodbreak}

\section{Model and results}\label{sec2}

\subsection{The model}
The model consists of $I$ customer classes and a single server.
A buffer with infinite room is dedicated to each customer class, and
upon arrival, customers are queued in the corresponding buffers.
Within each class, customers are served at the order of arrival.
The server may only serve the customer at the head of each line.
Moreover, processor sharing is allowed, and so the server
is capable of serving up to $I$ customers (of distinct classes) simultaneously.

The model is defined on a probability space $(\Om,\calF,\PP)$.
Expectation with
respect to $\PP$ is denoted by $\EE$.
The parameters and processes we introduce will depend on
an index $n\in\N$, that will serve as a scaling parameter.
Arrivals occur according to independent
renewal processes, and service times are independent and identically distributed
across each class.
Let $\I=\{1,2,\ldots, I\}$.
Let $\lambda^n_i >0, n\in\N, i\in\I$ be given parameters,
representing the \textit{reciprocal mean
inter-arrival times} of class-$i$ customers.
Given are $I$ independent sequence $\{\IA_i(l) \dvtx l\in\N\}_{i\in
\calI}$,
of positive i.i.d. random variables with mean $\E[\IA_i(1)]=1$ and
variance $\sig^2_{i,IA}=\operatorname{Var}(\IA_i(1))\in(0,\infty)$. With
$\sum_1^0=0$,
the number of arrivals of class-$i$ customers up to time $t$, for the
$n$th system,
is given by
\[
A^n_i(t)=\sup \Biggl\{ l\geq0 \dvtx \sum
_{k=1}^l\frac{\IA_i(k)}{\lambda
^n_i}\leq t \Biggr\},\qquad t\geq0.
\]
Similarly we consider another set of parameters $\mu^n_i>0, n\in\N,
i\in\calI$,
representing \textit{reciprocal mean service times}.
We are also given $I$ independent sequence $\{\ST_i(l) \dvtx l\in\N\}
_{i\in\calI}$ of positive i.i.d. random variables (independent also
of the sequences $\{\IA_i\}$)
with mean $\E[\ST_i(1)]=1$ and variance $\sig^2_{i,ST}=\operatorname
{Var}(\ST_i(1))\in(0,\infty)$.
The time required to complete the service of the $l$th class-$i$
customer is given
by $\ST_i(l)/\mu^n_i$, and the \textit{potential service time}
processes are defined as
\[
S^n_i(t)=\sup \Biggl\{ l\geq0 \dvtx \sum
_{k=1}^l\frac{\ST_i(k)}{\mu
^n_i}\leq t \Biggr\},\qquad t\geq0.
\]
We consider the \textit{moderate deviations rate parameters} $\{b_n\}
$, that form a sequence,
fixed throughout, with the property that $\lim b_n=\infty$ while $\lim
\frac{b_n}{\sqrt{n}}=0$, as $n\to\iy$.
The arrival and service parameters are assumed to satisfy the following
conditions.
As $n\to\iy$:
\begin{itemize}
\item$\frac{\lambda^n_i}{n}\to\lambda_i\in(0,\iy)$ and $\frac
{\mu^n_i}{n}\to\mu_i\in(0,\iy)$,
\item$\tilde\lambda^n_i:=\frac{1}{b_n\sqrt{n}}(\lambda
^n_i-n\lambda_i)\to\tilde\lambda_i\in(-\infty, \infty)$,
\item$\tilde\mu^n_i:=\frac{1}{b_n\sqrt{n}}(\mu^n_i-n\mu_i)\to
\tilde\mu_i\in(-\infty, \infty)$.
\end{itemize}
Also the system is assumed to be critically loaded in the sense that
$\sum_1^I \rho_i=1$ where $\rho_i=\frac{\lambda_i}{\mu_i}$ for
$i\in\calI$.

For $i\in\calI$, let $X^n_i$ be a process representing the number of
class-$i$ customers in the $n$th system.
With $\mathbb{S}=\{x=(x_1,\ldots,x_I)\in[0,1]^I\dvtx \sum x_i\le1\}$, let
$B^n$ be a process taking values in $\mathbb{S}$, whose $i$th
component represents
the fraction of effort devoted by the server to the class-$i$ customer
at the head of the line.
Then the number of service completions of class-$i$ jobs during
the time interval $[0,t]$ is given by
%
\begin{equation}
\label{01} D^n_i(t):=S^n_i
\bigl(T^n_i(t)\bigr),
\end{equation}
where
%
\begin{equation}
\label{02} T^n_i(t)=\int_0^t
B^n_i(s)\,ds
\end{equation}
is the time devoted to class-$i$ customers by time $t$.
The following equation follows from foregoing verbal description
%
\begin{equation}
\label{03} X^n_i(t)=X^n_i(0)+A^n_i(t)-S^n_i
\bigl(T^n_i(t)\bigr).
\end{equation}
For simplicity, the initial conditions $X^n_i(0)$ are assumed to be
deterministic.
Note that, by construction, the arrival and potential service processes
have RCLL paths,
and accordingly, so do $D^n$ and $X^n$.

The process $B^n$ is regarded as a control that is determined based on
observations
from the past (and present) events in the system. A precise definition
is as follows.
Fix $T>0$ throughout.
Given $n$, the process $B^n$ is said to be an \textit{admissible control}
if its sample paths lie in $\calD([0,T],\mathbb{S})$ and:
\begin{itemize}
\item it is adapted to the filtration
\[
\s\bigl\{A^n_i(s),S^n_i
\bigl(T^n_i(s)\bigr), i\in\I, s\leq t\bigr\},
\]
where $T^n$ is given by \eqref{02};
\item for $i\in\I$ and $t\ge0$, one has
%
\begin{equation}
X^n_i(t)=0 \mbox{ implies } B^n_i(t)=0,\label{04}
\end{equation}
where $X^n$ is given by \eqref{03}.
\end{itemize}
Denote the class of all admissible controls $B^n$ by $\Bf^n$.
Note that this class depends on $A^n$ and $S^n$, but we consider these processes
to be fixed. It is clear that this class is nonempty, as one may obtain
an admissible control,
for example, by setting $B^n=0$ identically.

We next introduce centered and scaled versions of the processes.
For $i\in\I$, let
%
\begin{eqnarray}
\label{22} \tilde{A}^n_i(t)&=&\frac{1}{b_n\sqrt{n}}
\bigl(A^n_i(t)-\lambda^n_i t
\bigr),\qquad \tilde{S}^n_i(t)=\frac{1}{b_n\sqrt{n}}
\bigl(S^n_i(t)-\mu^n_i t\bigr),
\nonumber
\\[-8pt]
\\[-8pt]
\nonumber
\tilde{X}^n_i(t)&=&\frac{1}{b_n\sqrt{n}}X^n_i(t).
\end{eqnarray}
It is easy to check from (\ref{03}) that
%
\begin{equation}
\tilde{X}^n_i(t)=\tilde{X}^n_i(0)
+ y^n_it +\tilde{A}^n_i(t)-
\tilde{S}^n_i\bigl(T^n_i(t)
\bigr) + Z^n_i(t), \label{eqn2}
\end{equation}
where we denote
%
\begin{equation}
\label{06} Z^n_i(t)=\frac{\mu^n_i}{n}
\frac{\sqrt{n}}{b_n}\bigl(\rho_i t-T^n_i(t)
\bigr),\qquad y^n_i=\tilde\lambda^n_i-
\rho_i\tilde\mu^n_i.
\end{equation}
Note that these processes have the property
%
\begin{equation}
\label{07} \sum_i\frac{n}{\mu^n_i}Z^n_i
\mbox{ starts from zero and is nondecreasing,}
\end{equation}
thanks to the fact that $\sum_iB^n_i\le1$ while $\sum_i\rho_i=1$.
Clearly $\tilde X^n_i$ is nonnegative, that is,
%
\begin{equation}
\label{08} \tilde X^n_i(t)\ge0,\qquad t\ge0, i\in\calI.
\end{equation}
We impose the following condition on the initial values:
\[
\tilde{X}^n(0)\to x\in\R^I_+ \qquad\mbox{as } n\to\infty.
\]

The scaled processes $(\tilde A^n,\tilde S^n)$ are assumed to satisfy a
\textit{moderate deviation
principle}. To express this assumption, let
$\Ir_k, k=1,2$, be functions on $\calD([0,T],\R^I)$ defined as
follows. For
$\psi=(\psi_1,\ldots,\psi_I)\in\calD([0,T],\R^I)$,
\[
\Ir_1(\psi)=\cases{ %
\displaystyle\frac{1}{2}\sum_{i=1}^I
\frac{1}{\lambda_i\sig^2_{i,IA}} \int_0^T\dot
\psi_i^2(s)\,ds, \vspace*{2pt}\cr
\qquad\hspace*{16pt} \mbox{if all } \psi_i
\mbox{ are absolutely continuous and } \psi(0)=0,
\vspace*{2pt}\cr
\infty, \qquad\mbox{otherwise},}
\]
and
\[
\Ir_2(\psi)=\cases{
\displaystyle\frac{1}{2}\sum_{i=1}^I
\frac{1}{\mu_i\sig^2_{i,ST}} \int_0^T\dot
\psi_i^2(s)\,ds,\vspace*{2pt}\cr
\qquad\hspace*{16pt} \mbox{if all } \psi_i
\mbox{ are absolutely continuous and } \psi(0)=0,
\vspace*{2pt}\cr
\infty, \qquad \mbox{otherwise}.}
\]
Let $\Ir(\psi)=\Ir_1(\psi^1)+\Ir_2(\psi^2)$
for $\psi=(\psi^1,\psi^2)\in\calD([0,T],\R^{2I})$. Note that $\Ir
$ is lower semicontinuous with compact level sets, properties used in
the sequel.

\begin{assumption}[(Moderate deviation principle)]\label{moderate}
The sequence
\[
\bigl(\tilde{A}^n,\tilde{S}^n\bigr)=\bigl(
\tilde{A}^n_1,\ldots, \tilde{A}^n_I,
\tilde{S}^n_1,\ldots,\tilde S^n_I
\bigr),
\]
satisfies the LDP with rate parameters $b_n$ and rate function
$\Ir$ in $\calD([0,T],\R^{2I})$; that is:
\begin{itemize}
\item for any closed set $F\subset\calD([0,T],\R^{2I})$
\[
\limsup\frac{1}{b_n^2}\log\p\bigl(\bigl(\tilde{A}^n,
\tilde{S}^n\bigr)\in F\bigr)\leq -\inf_{\psi\in F}\Ir(
\psi);
\]
\item for any open set $G\subset\calD([0,T],\R^{2I})$
\[
\liminf\frac{1}{b_n^2}\log\p\bigl(\bigl(\tilde{A}^n,
\tilde{S}^n\bigr)\in G\bigr)\geq -\inf_{\psi\in G}\Ir(
\psi).
\]
\end{itemize}
\end{assumption}

\begin{remark}[(Sufficient conditions)]\label{rem1}
It is shown in \cite{Puhal-Whitt} that each one of the following
statements is sufficient for
Assumption \ref{moderate} to hold:
\begin{itemize}
\item there exist constants $u_0>0$, $\beta\in(0,1]$ such that $\E
[e^{u_0(\IA_i)^\beta}], \E[e^{u_0(\ST_i)^\beta}]<\infty$, $i\in
\I$, and $b_n^{\beta-2}n^{\beta/2}\to\infty$;
\item for some $\eps>0$, $\E[(\IA_i)^{2+\eps}], \E[(\ST
_i)^{2+\eps}]<\infty$, $i\in\I$, and
$b_n^{-2}\log n\to\infty$.
\end{itemize}
\end{remark}

To present the risk-sensitive control problem, let $h$ and $g$ be
nonnegative, continuous functions from $\R^I_+$ to $\R$, monotone
nondecreasing with respect to the partial order
$a\le b$ if and only if $b-a\in\R^I_+$. Assume that $h, g$ have at
most linear growth, that is, there exist
constants $c_1, c_2$ such that
\[
g(x)+h(x)\leq c_1\|x\|+ c_2.
\]
Given $n$, the cost associated with the initial condition
$\tilde X^n(0)$ and control $B^n\in\Bf^n$ is given by
%
\begin{equation}
\label{05} J^n\bigl(\tilde{X}^n(0),B^n
\bigr)=\frac{1}{b_n^2} \log\E \bigl[e^{b^2_n[\int_0^Th(\tilde{X}^n(s))\,ds+g(\tilde
{X}^n(T))]} \bigr].
\end{equation}
The value function of interest is given by
\[
V^n\bigl(\tilde{X}^n(0)\bigr)=\inf_{B^n\in\Bf^n}
J^n\bigl(\tilde{X}^n(0),B^n\bigr).
\]

\subsection{A differential game}

We next develop a differential game for the limit behavior
of the above control problem.
Let $\theta=(\frac{1}{\mu^1},\ldots,\frac{1}{\mu^I})$ and
$y=(y_1,\ldots,y_I)$
where $y_i=\tilde\lambda_i-\rho_i\tilde\mu_i$.
Denote $\calP=\calC_0([0,T],\R^{2I})$,\vadjust{\goodbreak} the subset of $\calC
([0,T],\R^{2I})$ of functions starting from zero,
and
\[
\calE=\bigl\{\zeta\in\calC\bigl([0,T],\R^I\bigr)\dvtx \theta\cdot
\zeta\mbox{ starts from zero and is nondecreasing}\bigr\}.
\]
Endow both spaces with the uniform topology.
Let $\brho$ be the mapping from $\calD([0,T],\R^I)$ into itself
defined by
\[
\brho[\psi]_i(t)=\psi_i(\rho_it),\qquad t
\in[0,T], i\in\calI.
\]
Given $\psi=(\psi^1,\psi^2)\in\calP$ and $\zeta\in\calE$, the
\textit{dynamics associated
with initial condition $x$ and data $\psi,\zeta$} is given by
%
\begin{equation}
\varphi_i(t)= x_i+y_it+
\psi^1_i(t)-\brho\bigl[\psi^2
\bigr]_i(t)+\zeta_i(t),\qquad i\in\calI. \label{201}
\end{equation}
Note the analogy between the above equation and equation \eqref{eqn2},
and between the condition $\theta\cdot\zeta$ nondecreasing and
property \eqref{07}.
The following condition, analogous to property \eqref{08}, will also
be used, namely
%
\begin{equation}
\label{09} \varphi_i(t)\ge0, \qquad t\ge0, i\in\calI.
\end{equation}
The game is defined in the sense of Elliott and Kalton \cite{Elli-Kal},
for which we need the notion of strategies.
A measurable mapping $\alpha\dvtx \calP\to\calE$ is called a \textit
{strategy for the minimizing player}
if it satisfies a causality property. Namely,
for every $\psi=(\psi^1,\psi^2), \tilde\psi=(\tilde\psi^1,\tilde
\psi^2) \in\calP$
and $t\in[0,T]$,
%
\begin{eqnarray}
\label{11} \bigl(\psi^1,\brho\bigl[\psi^2
\bigr]\bigr) (s)=\bigl(\tilde\psi^1,\brho\bigl[\tilde
\psi^2\bigr]\bigr) (s)
\nonumber
\\[-8pt]
\\[-8pt]
 \eqntext{\mbox{for all } s\in[0,t]
\mbox{ implies } \alpha[\psi](s)=
\alpha[\tilde\psi](s) \mbox{ for all }  s\in[0,t].}
\end{eqnarray}
Given an initial condition $x$,
a strategy $\al$ is said to be \textit{admissible} if, whenever $\psi
\in\calP$ and
$\zeta=\al[\psi]$, the corresponding dynamics \eqref{201} satisfies
the nonnegativity
constraint \eqref{09}. The set of all admissible strategies for the minimizing
player is denoted by $A$ (or, when the dependence on the initial
condition is
important, $A_x$).
Given $x$ and $(\psi,\zeta)\in\calP\times\calE$, we define the
cost by
\[
c(\psi,\zeta)=\int_0^Th\bigl(\varphi(t)
\bigr)\,dt+g\bigl(\varphi(T)\bigr)-\Ir(\psi),
\]
where $\varphi$ is the corresponding dynamics.
The value of the game is defined by
\[
V(x)=\inf_{\alpha\in A_x}\sup_{\psi\in\calP}c\bigl(\psi,
\alpha[\psi]\bigr).
\]

\subsection{Main result}

For $w\in\R_+$, denote
%
\begin{eqnarray}
\label{84} h^*(w)&=&\inf\bigl\{h(x) \dvtx x\in\R^I_+, \theta\cdot
x=w\bigr\},
\nonumber
\\[-8pt]
\\[-8pt]
\nonumber
 g^*(w)&=&\inf\bigl\{g(x) \dvtx x\in\R^I_+, \theta\cdot x=w
\bigr\}.
\end{eqnarray}\eject\noindent
We need the following assumption. It is similar to the one imposed in
\cite{atasol,atagur},
where an analogous many-server model is treated in a diffusion regime.
%
\begin{assumption}[(Existence of a continuous minimizing curve)]\label{mini}
There exists a
continuous map $f \dvtx \R_+\to\R^I_+$ such that for all $w\in\R_+$,
\[
\theta\cdot f(w)=w,\qquad h^*(w)=h\bigl(f(w)\bigr),\qquad g^*(w)=g\bigl(f(w)\bigr).
\]
\end{assumption}
As far as solving the game is concerned, this assumption is not
required at all; see Remark \ref{rem2}.
It is important in the proof of asymptotic optimality.
The fact that the same function $f$ serves as a minimizer for both $h$
and $g$
may seem to be too strong.
We comment in Remark \ref{rem3} on what is involved in relaxing this
assumption.
%
\begin{example}\label{mini-rem}
(a) The linear case: $h(x)=\sum c_ix_i$ and $g(x)=\sum d_ix_i$, for some
nonnegative constants $c_i$, $d_i$. If we require
that $c_I\mu_I=\min_i c_i\mu_i$ and $d_I\mu_I=\min_i d_i\mu_i$, then
the assumption holds with $f(w)=(0,\ldots,0, w\mu_I)$.
This is the case considered in Section~\ref{sec5}.

(b) If $h$ is nondecreasing, homogeneous of degree $\alpha,
0<\alpha\leq1$, and $x^*\in\operatorname{argmin}\{h(x) \dvtx \theta\cdot x=1\}
$, it
is easy to check that $f(w)=wx^*$ satisfies the above assumption
provided $g=dh$ for some nonnegative constant $d$.
\end{example}

\begin{assumption}[(Exponential moments)]\label{unbounded}
For any constant $K$,
\[
\limsup_{n\to\iy}\frac{1}{b^2_n}\log\E\bigl[e^{b^2_nK(\|\tilde{A}^n\|
^*_T+\|\tilde{S}^n\|^*_T)}
\bigr] <\infty.
\]
\end{assumption}
A sufficient condition for the above is as follows (see the \hyperref[app]{Appendix}
for a proof).

\begin{proposition}\label{prop-bd}
If there exists $u_0>0$ such that $\E[e^{u_0 \IA_i}]$ and $\E[e^{u_0
\ST_i}]$, $i\in\I$, are finite,
then Assumption \ref{unbounded} holds.
\end{proposition}
Note that taking $\beta=1$ in Remark \ref{rem1} shows that the
hypothesis of Proposition~\ref{prop-bd}
is sufficient for Assumption \ref{moderate} as well.

Our main result is the following:
%
\begin{theorem}\label{main}
Let Assumptions \ref{moderate} and \ref{mini} hold. If either $g$ or
$h$ is unbounded,
let also Assumption \ref{unbounded} hold. Then
$\lim_{n\to\iy} V^n(\tilde{X}^n(0))=V(x)$.
\end{theorem}

\begin{remark}[(An equivalent game)]
There is a simpler, equivalent formulation of the game, which avoids
the use of the time scaling operator $\brho$
(both formulations will be used in the proofs).
Define a functional $\bar\Ir(\psi)=\bar\Ir_1(\psi^1)+\bar\Ir
_2(\psi^2)$ on $\calD([0,T], \R^{2I})$,
where $\bar\Ir_k, k=1,2$, are functionals
on $\calD([0,T],\R^I)$ given by $\bar\Ir_1=\Ir_1$, and,
for $\psi=(\psi_1,\ldots,\psi_I)\in\calD([0,T],\R^I)$,
\[
\bar\Ir_2(\psi)=\cases{ %
\displaystyle\frac{1}{2}\sum_{i=1}^I
\frac{1}{\rho_i\mu_i\sig^2_{i,ST}} \int_0^T\dot
\psi_i^2(s)\,ds, \vspace*{2pt}\cr
\qquad\hspace*{16pt}\mbox{if all } \psi_i
\mbox{ are absolutely continuous and } \psi(0)=0,
\vspace*{2pt}\cr
\infty,\qquad \mbox{otherwise}.}
\]
The dynamics of the game $\bar\varphi$ are now
\[
\bar{\varphi}_i(t)=x_i=y_it+
\psi^1(t)-\psi^2(t)+\zeta_i(t)\geq0.
\]
A strategy $\al$ should now satisfy the following version of the
causality property:
\[
\mbox{$\psi(s)=\psi(s)$ for all $s\in[0,t]$ implies $\alpha[\psi](s)=\alpha[
\tilde\psi](s)$ for all $s\in[0,t]$.}
\]
Denote the set of all such strategies by $\bar A_x$.
Given $x$ and $(\psi,\zeta)\in\calP\times\calE$, let
\[
\bar c(\psi,\zeta)=\int_0^Th\bigl(\bar
\varphi(t)\bigr)\,dt+g\bigl(\bar\varphi (T)\bigr)-\bar\Ir(\psi),
\]
where $\bar\varphi$ is as above.
Then the value of the game can also be defined as
\[
V(x)=\inf_{\alpha\in\bar A_x}\sup_{\psi\in\calP}\bar c\bigl(\psi,
\alpha[\psi]\bigr).
\]
\end{remark}

\begin{remark}[(Possible extensions)] Our main results can be extended in
various ways.
The following two are relatively easy. We do not provide the proofs because
we aim at keeping these aspects as simple as possible in this paper.

(a) The moderate deviation principle (Assumption \ref{moderate}),
that forms the basis of the asymptotic analysis,
is proved in \cite{Puhal-Whitt} to hold for a sequence of renewal
processes in a more general formulation, namely that of a triangular array.
Our results can be extended to cover this formulation.

(b) The assumption that the $2I$ service and arrival
processes are mutually
independent leads to the form $\mathbb{I}(\psi)=\int_0^TF[\dot\psi
(s)]\,ds$ [if $\psi$ is absolutely continuous and $\psi(0)=0$; $\iy$
otherwise] of the rate function,
where $F$ is a weighted sum of squares. Our results can be extended to cover
dependence structure such as where $F$ is a positive definite quadratic form.
\end{remark}

\section{Solution of the game}\label{sec3}

In this section we find a minimizing strategy for $V$,
under Assumption \ref{mini}, following an idea from \cite{Har-Van}.
Throughout this section, the initial condition $x$ is fixed.
Consider the one-dimensional \textit{Skorohod map} $\Gam$ from $\calD
([0, T], \R)$ to itself given by
%
\begin{equation}
\Gam[z](t)=z(t)-\inf_{s\in[0,t]}\bigl[z(s)\wedge0\bigr],\qquad t\in
[0,T].\label{skoro}
\end{equation}
Clearly, $\Gam[z](t)\ge0$ for all $t$. Let also
\[
\bar\Gam[z](t)=-\inf_{s\in[0,t]}\bigl[z(s)\wedge0\bigr],\qquad t\in[0,T].
\]
It is clear from the definition that, for $z,w\in\calD([0,T], \R)$
%
\begin{equation}
\sup_{[0,T]}\bigl|\Gam[z]-\Gam[w]\bigr|\leq2 \sup_{[0,T]}|z-w|.\label{302}
\end{equation}
The construction below is based on the mapping $\Gam$ and the function $f$
from Assumption \ref{mini}.
Recall from \eqref{201} that for $\psi=(\psi^1,\psi^2)\in\calP$
and $\zeta\in\calE$,
the dynamics of the differential game is given by
\[
\varphi=\xi+\zeta,
\]
where
\[
\xi(t)=x+yt+\psi^1(t)-\brho\bigl[\psi^2\bigr](t),\qquad t
\in[0,T].
\]
We associate with each $\psi\in\calP$ a 4-tuple
$(\bph[\psi],\bxi[\psi],\bzeta[\psi],\bw[\psi])$ given by
%
\begin{eqnarray}
\label{70} \bxi[\psi](t) &=& x+yt+\psi^1(t)-\brho\bigl[
\psi^2\bigr](t),\qquad t\in[0,T],
\\
\label{71} \bw[\psi] &=&\Gam\bigl[\theta\cdot\bxi[\psi]\bigr],
\\
\label{72} \bph[\psi] &=& f\bigl(\bw[\psi]\bigr),
\\
\label{73} \bzeta[\psi] &=& \bph[\psi]-\bxi[\psi].
\end{eqnarray}
Sometimes we also use the notation
%
\begin{eqnarray}
\label{80} \hat{\bxi}[\psi](t) &=& x+yt+\psi^1(t)-
\psi^2(t),\qquad t\in[0,T],
\\
\label{81} \hat{\bw}[\psi] &= &\Gam\bigl[\theta\cdot\hat{\bxi}[\psi]\bigr],
\\
\label{82} \hat{\bph}[\psi] &=& f\bigl(\hat{\bw}[\psi]\bigr),
\\
\label{83} \hat{\bzeta}[\psi] &= &\hat{\bph}[\psi]-\hat{\bxi}[\psi].
\end{eqnarray}
Note that $\bzeta[\psi^1,\psi^2]=\hat{\bzeta}[\psi^1,\brho[\psi^2]]$.

As we state in the result below, $\bzeta$ is an optimal strategy.
Now, the state variable $\varphi$ generally lies in $I$ dimensions.
But under the solution provided by $\bzeta$, namely when
$(\varphi,\xi,\zeta,w)=(\bph,\bxi,\bzeta,\bw)[\psi]$, one has
$\varphi=f(w)$,
and so the state variable lies on a \textit{one-dimensional} manifold,
and is dictated solely by
the one-dimensional object $w$, that represents workload.
This dimensionality reduction owes to
the fact that, in the scaling limit, a proper allocation of effort at
the server
can drive the state variable $\varphi$ instantaneously to the location
$\varphi=f(w)$.
As far as the literature on heavy traffic limits at the \textit
{diffusion scale}
is concerned, the instantaneous mobility as well as the reduction to a
problem based on
the workload dimension (called workload reduction)
are well known for this and related models. See, for example, the
explanation of a similar phenomenon in \cite{man-sto}, and
general\vadjust{\goodbreak}
results on
workload reduction in \cite{Har-Wil}.
Our results thus establish the validity of workload reduction at the
\textit{MD scale},
for the model under study.

\begin{proposition}
\label{prop1}
Let Assumption \ref{mini} hold. Then $\bzeta$ is an admissible
strategy. Moreover,
it is a minimizing strategy, namely
%
\begin{equation}
V(x)=\sup_{\psi\in\calP} c\bigl(\psi, \bzeta[\psi]\bigr).
\label{305}
\end{equation}
\end{proposition}
\begin{pf}
Let us show that $\bzeta$ is an admissible strategy. Let $\psi\in
\calP$ be given
and denote $(\varphi,\xi,\zeta,w)=(\bph,\bxi,\bzeta,\bw)[\psi]$.
Then $\varphi=\xi+\zeta$, and multiplying \eqref{73} by $\theta$,
\[
\theta\cdot\zeta=w-\theta\cdot\xi=\bar\Gam[\theta\cdot\xi].
\]
Since $\theta\cdot\xi(0)=\theta\cdot x\ge0$, it follows that
$\theta\cdot\zeta(0)=0$.
Moreover, by definition of $\bar\Gam$, $\theta\cdot\zeta$ is
nondecreasing. This shows
$\zeta\in\calE$. The causality property \eqref{11} follows directly
from an analogous property of
$\bar\Gam$. Next, $w(t)\ge0$ for all $t$, and, by definition, $f$ maps
$\R_+$ to $\R_+^I$, whence $\varphi(t)\in\R_+^I$ for all $t$. This
shows that
$\bzeta$ is an admissible strategy.

Now we check that $\bzeta$ is indeed a minimizing strategy.
This is based on the minimality property of the Skorohod map; see, for
example, \cite{cheman}, Section~2.
Namely, if $z,r\in\calD([0,T]\dvtx \R)$, $r$ is nonnegative and nondecreasing,
and $z(t)+r(t)\ge0$ for all $t$, then
\[
z(t)+r(t)\ge\Gam[z](t),\qquad t\in[0,T].
\]
Let $\alpha\in A$ be any admissible strategy. Given $\psi$, let
$(\varphi,\xi,\zeta,w)$ be as before.
The dynamics corresponding to $\psi$ and $\tilde\zeta:=\al[\psi]$
is given
by $\tilde\varphi=\xi+\tilde\zeta$. Since $\al$ is an admissible
strategy, we have that
\[
\theta\cdot\tilde\varphi=\theta\cdot\xi+\theta\cdot\tilde \zeta\ge0,
\]
and $\theta\cdot\tilde\zeta$ is
nonnegative and nondecreasing. Thus by the above minimality property,
\[
\theta\cdot\tilde\varphi(t)\ge\Gam[\theta\cdot\xi](t)=w(t),\qquad t\in[0,T].
\]
By monotonicity of $h$, \eqref{72} and Assumption \ref{mini},
\begin{eqnarray}
\label{16}
h\bigl(\tilde\varphi(t)\bigr)& \geq&\inf\bigl\{h(q) \dvtx
\theta\cdot q=\theta \cdot\tilde\varphi(t)\bigr\}
\nonumber
\\[-8pt]
\\[-8pt]
\nonumber
&\geq&\inf\bigl\{h(q) \dvtx
\theta\cdot q=w(t)\bigr\} = h\bigl(f\bigl(w(t)\bigr)\bigr)=h\bigl(\varphi(t)
\bigr).
\end{eqnarray}
A similar estimate holds for $g$, namely
%
\begin{equation}
\label{17} g\bigl(\tilde\varphi(T)\bigr)\ge g\bigl(\varphi(T)\bigr).
\end{equation}
As a result,
\[
\sup_{\psi\in\calP}c\bigl(\psi, \alpha[\psi]\bigr)\geq \sup
_{\psi\in\calP}c\bigl(\psi,\bzeta[\psi]\bigr).
\]
This proves that $\bzeta$ is a minimizing strategy; namely \eqref
{305} holds.
\end{pf}

\begin{remark}[(Beyond Assumption \ref{mini})]
\label{rem2}
(a) The game can be solved without Assumption \ref{mini}.
Owing to the continuity of $h$ and $g$ and using a measurable selection result
such as Corollary 10.3 in the Appendix of \cite{ethkur}, there exist measurable
functions $f_h$ and $f_g$ mapping $\R_+$ to $\R_+^I$ such that for
all $w\in\R_+$,
%
\begin{eqnarray}
\label{85} \theta\cdot f_h(w)&=&\theta\cdot f_g(w)=w,\qquad
h^*(w)=h\bigl(f_h(w)\bigr),
\nonumber
\\[-8pt]
\\[-8pt]
\nonumber
 g^*(w)&=&g\bigl(f_g(w)\bigr),
\end{eqnarray}
where we recall the definition \eqref{84} of $h^*$ and $g^*$.
To construct a minimizing strategy, let $\bxi$ and $\bw$ be as in
\eqref{70}--\eqref{71}. Instead of \eqref{72}, consider
%
\begin{equation}
\label{87} \bph[\psi](t)=\cases{ f_h\bigl(\bw[\psi](t)\bigr), &\quad
$t\in[0,T)$,\vspace*{2pt}
\cr
f_g\bigl(\bw[\psi](T)\bigr), & \quad $t=T$.}
\end{equation}
Then define $\bzeta$ as in \eqref{73} accordingly ($\calE$ and
$\calP$ will also
change accordingly).
The proof of Proposition \ref{prop1} goes through with almost no change.
Indeed, the continuity of $f$ is not used in this proof, and inequalities
\eqref{16} and \eqref{17} can be obtained by working with $f_h$ and $f_g$,
respectively, instead of $f$.

(b) Although the continuity that is a part of in
Assumption \ref{mini} is irrelevant
for the game, it will be used in the convergence argument leading to
the asymptotic optimality result (Theorem \ref{th42}). One may,
however, consider
a relaxation of Assumption \ref{mini} as follows: There exist
continuous functions
$f_h$ and $f_g$ satisfying~\eqref{85} above. Under this relaxed assumption,
given a continuous path $\psi\in\calP$,
the corresponding dynamics $\varphi=\bph[\psi]$, with $\bph$ as in
\eqref{87},
may then have a jump at time $T$. The jump makes it more complicated
to obtain convergence in Theorem~\ref{th42}. We discuss this issue in
Remark \ref{rem3}.
\end{remark}

\textit{Extension and some properties of $\bzeta$}.
As a strategy, $\bzeta$ is defined on $\calP$
[recall $\calP=\calC_0([0, T], \R^{2I})$]. We extend $\hat{\bzeta
}$ and $\bzeta$ to
\[
\bar\calP=\calD\bigl([0,T],\R^{2I}\bigr),
\]
using the same definitions \eqref{73} and \eqref{83}.
Some useful properties related to this map are stated in the following result.
Given a map $m\dvtx [0, T]\to\R^k$, some $k\in\N$, and
$\eta>0$, define the $\eta$-oscillation of $m$ as
\[
\osc_\eta(m)=\sup\bigl\{\bigl\|m(s)-m(t)\bigr\| \dvtx |s-t|\leq\eta, s,t\in[0,
T]\bigr\}.
\]
For $\kappa>0$, define (with $\|\cdot\|^*=\|\cdot\|^*_T$)
%
\begin{equation}
\label{19} \calD(\kappa)=\bigl\{\psi=\bigl(\psi^1,
\psi^2\bigr)\in\bar\calP \dvtx \bigl\|\psi^1\bigr\|^*+\bigl\|
\psi^2\bigr\|^*\leq\kappa \mbox{ and } \bxi[\psi ](0)\in\R^I_+
\bigr\}.
\end{equation}

\begin{proposition}
\label{prop2}
Let Assumption \ref{mini} hold.

\begin{longlist}[(iii)]
\item[(i)] Given $\xi,\zeta\in\calD([0,T],\R^I)$,
$\varphi(t)=\xi(t)+\zeta(t)\in\R_+^I$,
$\theta\cdot\zeta$ nonnegative and nondecreasing, one has
%
\begin{equation}
\label{18} j\bigl(\varphi(t)\bigr)\ge j\bigl(f\bigl(\Gam[\theta\cdot\xi](t)
\bigr)\bigr)\qquad \mbox{for } j=h,g.\vadjust{\goodbreak}
\end{equation}
\item[(ii)] There exist constants $\gamma_0$ and $\gamma_1$ such that for
$\psi\in\bar\calP$,
%
\begin{equation}
\bigl\|\hat{\bzeta}[\psi](t)\bigr\|\leq\gamma_0\bigl(\bigl\|\psi^1
\bigr\|^*_t+\bigl\|\psi^2\bigr\| ^*_t\bigr)+
\gamma_1. \label{311}
\end{equation}
\item[(iii)] For $\psi, \tilde\psi\in\calD(\kappa)$, given $\eps>0$
there exists $\delta_1>0$ such that
%
\begin{equation}
\bigl\|\hat{\bzeta}[\psi] -\hat{\bzeta}[\tilde\psi]\bigr\|^*\le\eps \qquad\mbox{provided } \bigl\|
\psi ^1-\tilde\psi^1\bigr\|^* +\bigl\|\psi^2-\tilde
\psi^2\bigr\|^*\leq\delta_1. \label{315}
\end{equation}

\item[(iv)]
For any $\psi\in\calD(\kappa)$, given $\eps>0$ there exist
$\delta>0$ and $\eta>0$ such that
%
\begin{equation}
\osc_\eta\bigl(\hat{\bzeta}[\psi]\bigr)\le\eps\qquad \mbox{provided } \osc
_\eta(\psi)\le\delta. \label{316}
\end{equation}
\end{longlist}
\end{proposition}
\begin{pf}
(i) The argument leading to \eqref{16} and \eqref{17} is seen to be applicable
for this extended map, giving \eqref{18}.\vspace*{-6pt}
\begin{longlist}[(iii)]
\item[(ii)] Denote $\theta_*=\min_{i\in\I}\theta_i$ and $\theta^*=\max_{i\in\calI}\theta_i$.
Then Assumption \ref{mini}
implies that $\|f(w)\|\leq\frac{1}{\theta_*}w$ for $w\geq0$. Let
$\gamma_0=\sqrt{I}(\frac{2\theta^*}{\theta_*}+1)$ and
$\gamma_1=\gamma_0\sum_{i=1}^I (x_i+T|y_i|)$.
Then for $t\in[0,T]$, using \eqref{80}--\eqref{83}, \eqref{311} holds.

\item[(iii)]
Using (\ref{82}) and (\ref{311}), for every $\kappa$ there exists a constant
$\beta=\beta(\kappa)$ such that, for all $\psi\in\calD(\kappa)$,
\begin{eqnarray*}
\bigl\|\hat{\bzeta}[\psi]\bigr\|^* &\leq& \beta(\kappa),
\\
\bigl|\hat\bw[\psi]\bigr|^* &\leq&\beta(\kappa).
\end{eqnarray*}
Thus given $\eps>0$ we can find $\delta=\delta(\kappa,\eps)$ such that
$
\|f(w_1)-f(w_2)\|<\frac{\eps}{2}  \mbox{ if } |w_1-w_2|\leq\delta \mbox{ and }
w_i\in[0, \beta(\kappa)]$.
Also using the relation $\hat{\bw}[\psi]=\Gam[\theta\cdot\hat
{\bxi}[\psi]]$
and the Lipschitz property of $\hat{\bxi}$,
we have for $\psi,\tilde\psi\in\bar\calP$
\[
\bigl|\hat{\bw}[\psi]-\hat{\bw}[\tilde\psi]\bigr|^*\leq c_1\bigl(\bigl\|\psi
^1-\tilde\psi^1\bigr\|^*+\bigl\|\psi^2-\tilde
\psi^2\bigr\|^*\bigr)
\]
for some constant $c_1$.
Choosing $\delta_1=\delta_1(\kappa,\eps)$ sufficiently small,
for $\psi, \tilde\psi\in\calD(\kappa)$ we have,
with $\varphi$ and $\tilde\varphi$ denoting the dynamics
corresponding to
$(\psi,\hat{\bzeta}[\psi])$ and, respectively, $(\tilde\psi,\hat
{\bzeta}[\tilde\psi])$,
\[
\|\varphi-\tilde\varphi\|^*\leq\frac{\eps}{2} \qquad\mbox{if } \bigl\|
\psi^1-\tilde\psi^1\bigr\|^*+\bigl\|\psi^2-\tilde
\psi^2\bigr\|^*\leq\delta_1.
\]
Using the above estimate and (\ref{83}) we have \eqref{315}.

\item[(iv)]
Property \eqref{316} follows directly from the definition of $\Gam$,
definitions \eqref{80}--\eqref{83} and the continuity of $f$.\quad\qed
\end{longlist}
\noqed\end{pf}

\section{Proof of Theorem \texorpdfstring{\protect\ref{main}}{2.1}}\label{sec4}

\subsection{Lower bound}\label{sec41}

\begin{theorem}\label{th41}
$\!\!\!$Let Assumptions \ref{moderate} and \ref{mini} hold. Then $\liminf
V^n(\tilde{X}^n(0))\geq V(x)$.
\end{theorem}

In the proof, we choose any path $\tilde\psi\in\calP$ and show that
for any nearly optimal policy, the paths $\tilde X^n(\cdot)$ can be
controlled suitably for $(\tilde A^n, \tilde S^n)$ close to $\tilde
\psi$. We\vadjust{\goodbreak} find a constant $G>0$ such that for $\theta\cdot Z^n>G$ the
lower bound becomes trivial by using the
monotonicity of $h$ and $g$, and for $\theta\cdot Z^n\leq G$, the
optimality of $\bzeta$ gives the required estimates.

\begin{pf*}{Proof of Theorem \ref{th41}}
Fix $\tilde\psi=(\tilde\psi^1,\tilde\psi
^2)\in\calP$.
Let $d(\cdot,\cdot)$ be as in \eqref{74}.
Define, for $r>0$,
\[
\A_r=\bigl\{\psi\in\calD\bigl([0,T], \R^{2I}\bigr) \dvtx
d(\psi,\tilde\psi)<r\bigr\}.
\]
Since $\tilde\psi$ is continuous, for any $r_1\in(0, 1)$ there
exists $r>0$ such that
%
\begin{equation}
\psi\in\A_r \mbox{ implies } \|\psi-\tilde\psi\|^*<r_1.
\label{411}
\end{equation}
Define $\theta^n=(\frac{n}{\mu^n_1},\frac{n}{\mu^n_2},\ldots
,\frac{n}{\mu^n_I})$. Then
$\theta^n\to\theta$ as $n\to\infty$.
Now, given $0<\eps<1$, choose a sequence of policies $\{B^n\}$ such that
%
\begin{equation}
\label{75} V^n\bigl(\tilde{X}^n(0)\bigr) + \eps>
J^n\bigl(\tilde{X}^n(0), B^n\bigr) \quad\mbox{and}\quad
B^n\in\Bf^n \qquad\mbox{for all } n.
\end{equation}
Recall that
%
\begin{equation}
J^n\bigl(\tilde{X}^n(0), B^n\bigr)=
\frac{1}{b_n^2}\log\E\bigl[e^{b_n^2[\int
_0^Th(\tilde{X}^n(s))\,ds
+g(\tilde{X}^n(T))]}\bigr],\label{412}
\end{equation}
where
%
\begin{eqnarray}
\tilde{X}^n_i(t)& =& \tilde{X}^n_i(0)+y^n_it+
\tilde{A}^n_i(t)-\tilde{S}^n_i
\bigl(T^n_i(t)\bigr) +Z^n_i(t),
\label{413}
\\
\label{53} Z^n_i(t)&=&\frac{\mu^n_i}{n}
\frac{\sqrt{n}}{b_n}\bigl(\rho_i t-T^n_i(t)
\bigr), \qquad T^n_i(t) = \int_0^tB^n_i(s)\,ds.
\end{eqnarray}
For $G>0$, define a random variable $\tau_n$ by
\begin{eqnarray*}
\tau_n&=&\inf\bigl\{t\ge0\dvtx \theta^n\cdot
Z^n(t)>G\bigr\}\w T \\
&\equiv&\inf \Biggl\{t\geq0 \dvtx
\frac{\sqrt{n}}{b_n} \Biggl(t-\sum_{i=1}^IT^n_i(t)
\Biggr)>G \Biggr\} \wedge T.
\end{eqnarray*}
By \eqref{07}, $\theta^n\cdot Z^n$
is nondecreasing and continuous and hence
\begin{eqnarray*}
\theta^n\cdot Z^n(t) &\leq& G \qquad\mbox{for } t\leq
\tau_n,
\\
\theta^n\cdot Z^n(t) &> & G\qquad \mbox{for } t>
\tau_n.
\end{eqnarray*}
Consider the event $(\tilde{A}^n,\tilde{S}^n)\in\A_r$. Under this event,
for $t>\tau_n$,
\[
\theta^n\cdot\tilde{X}^n(t) \geq-\bigl\|\theta^n
\bigr\|\bigl(\kappa_0+2\|\tilde\psi\|^*\bigr)+G,
\]
where $\kappa_0$ is a constant (not depending on $n$ or $G$), and
we used (\ref{411}) and the boundedness of $\tilde X^n(0)$ and
$\tilde\lambda^n_i-\rho_i\tilde\mu^n_i$.
Since also $\theta^n$ converges,
we can choose a constant
$\kappa_1$
such that, on the indicated event,
%
\begin{equation}
\theta^n\cdot\tilde{X}^n(t)\geq-\kappa_1+G,\qquad
t>\tau_n. \label{414}
\end{equation}

Next, let $w=\bw[\tilde\psi]$,
$\varphi=\bph[\tilde\psi]$, $\zeta=\bzeta[\tilde\psi]$; see
\eqref{70}--\eqref{73}.
Note that $\varphi$ is the dynamics corresponding to $(\tilde\psi
,\zeta)$, namely
%
\begin{equation}
\label{12} \varphi_i(t)=x_i+y_it+
\tilde\psi^1_i(t)-\tilde\psi^2_i(
\rho_it) +\zeta_i(t).\vadjust{\goodbreak}
\end{equation}
For any $\kappa>0$ define a compact set $Q(\kappa)$ as
\[
Q(\kappa)=\bigl\{q\in\R^I_+ \dvtx 2q\cdot\theta\leq\kappa\bigr\}.
\]
Choose $\kappa$ large enough so that
\[
h(z)\geq\bigl|h\bigl(\varphi(\cdot)\bigr)\bigr|^*_T \quad\mbox{and}\quad g(z)\geq g
\bigl(\varphi(T)\bigr)
\]
for all $z\in Q^c(\kappa)$. To see that this is possible
note that $h(f(\cdot))$ is nondecreasing, and for $z\in Q^c(\kappa)$
\[
h(z)\geq\min\bigl\{h(q) \dvtx \theta\cdot q=\theta\cdot z\bigr\} =h\bigl(f(
\theta\cdot z)\bigr),
\]
where we use the definition of $f$. Thus
\[
h(z)\ge h\bigl(f(\kappa/2)\bigr),
\]
where we use the monotonicity of $h(f(\cdot))$.
Since $\tilde\psi(t)$, $t\in[0,T]$, is bounded, so is $w(t)$, $t\in
[0,T]$, by continuity of
$\Gam$. Choosing $\kappa=2|w|^*_T$ and using again the monotonicity
of $h(f(\cdot))$,
gives the claimed inequality for $h$. A similar argument applies for $g$.

Since $\theta_*:=\min_i\theta_i>0$,
we can choose $n_0$ large enough to ensure that $\theta^n_i
\leq2\theta_i$ for all $i\in\I$ and $n\geq n_0$.
Now if we choose $G$ in (\ref{414}) large enough so that $-\kappa
_1+G>\kappa$,
we have for $t>\tau_n, n\geq n_0$,
\[
2\theta\cdot\tilde{X}^n(t)\geq\theta^n\cdot
\tilde{X}^n(t)>\kappa,
\]
and hence by our choice of $\kappa$ we have on the indicated event,
for all $t>\tau_n$,
%
\begin{eqnarray}\label{444}
h\bigl(\tilde{X}^n(t)\bigr)\geq\bigl|h(\varphi)\bigr|^* \quad\mbox{and}\quad g\bigl(
\tilde{X}^n(t)\bigr) \geq g\bigl(\varphi(T)\bigr)
\nonumber
\\[-8pt]
\\[-8pt]
\eqntext{\mbox{for all
sufficiently large } n.}
\end{eqnarray}

Now we fix $G$ as above and consider $t\leq\tau_n$, on the same event
$(\tilde{A}^n,\tilde{S}^n)\in\A_r$.
The nonnegativity of $\tilde X^n_i$ and \eqref{413} imply a lower
bound on each of
the terms $Z^n_i$, namely
\[
Z^n_i(t)\geq-\tilde{X}^n(0)-y^n_it
-\tilde{A}^n_i(t)+\tilde{S}^n_i
\bigl(T^n_i(t)\bigr).
\]
Therefore using (\ref{411}) there exists a constant $\kappa_2$
such that for all sufficiently large $n$,
$
Z^n_i(t)\geq-\kappa_2$.
Combining this with the definition of $\tau_n$ in terms of $G$, we have
for $t\leq\tau_n$ and all large $n$,
%
\begin{equation}
\bigl\|Z^n(t)\bigr\|\leq\kappa_3. \label{417}
\end{equation}
Consider the stochastic processes $\Ps^n, \tilde\Ps^n, \tilde Z^n$,
with values in $\R^I$,
\begin{eqnarray*}
\Ps^n_i(t) &=& \tilde{A}^n_i(t
\w\tau_n),
\\
\tilde{\Ps}^n_i(t) &=& x_i-
\tilde{X}^n_i(0)+\bigl(y_i-y^n_i
\bigr)t+\tilde {S}^n_i\bigl(T^n_i(t
\w\tau_n)\bigr) -\bigl(1-\mu_i\theta^n_i
\bigr)Z^n_i(t\w\tau_n),
\\
\tilde Z^n_i(t) &=& \mu_i
\theta^n_iZ^n_i(t).
\end{eqnarray*}
Then by (\ref{413}),
%
\begin{equation}
\tilde{X}^n_i(t) = x_i+y_it+
\Ps^n_i(t)-\tilde\Ps^n_i(t)+
\tilde Z^n_i(t),\qquad t\in[0,\tau_n].
\label{418}
\end{equation}
Note that $\Ps^n, \tilde{\Ps}^n$ have RCLL sample paths, and consider
$\Ph^n=\hat{\bph}[\Ps^n,\tilde\Ps^n]$. Then
%
\begin{equation}
\label{13} \Ph^n(t)=x+yt+\Ps^n(t)-\tilde{
\Ps}^n(t)+\hat{\bzeta}\bigl[\Ps ^n,\tilde\Ps^n
\bigr](t).
\end{equation}
Let us now apply Proposition \ref{prop2}(i)
with $\xi(t)=x+yt+\Ps^n(t)-\tilde\Ps^n(t)$ and $\zeta=\tilde Z^n$.
Note that $\tilde X^n=\xi+\zeta$ takes values in $\R_+^I$, by
definition, and that
$\theta\cdot\tilde Z^n$ is nonnegative and nondecreasing, by \eqref{07}.
Moreover, by definition of $\hat{\bph}$ [see \eqref{80}--\eqref
{82}], $\Ph^n=f(\Gam[\theta\cdot\xi])$.
Hence \eqref{18} gives
%
\begin{equation}
\label{14} h\bigl(\tilde{X}^n(t)\bigr)\geq h\bigl(
\Ph^n(t)\bigr) \quad\mbox{and}\quad g\bigl(\tilde {X}^n(t)\bigr)\geq
g\bigl(\Ph^n(t)\bigr),\qquad t\in[0,\tau_n].
\end{equation}

Let $\kappa_4=\|\tilde\psi\|^*$. By (\ref{411}), on the indicated event,
$(\tilde{A}^n, \tilde{S}^n)\in\calD(2(1+\kappa_4))$ where we
recall definition \eqref{19}.
Note that $x+\Ps^n(0)-\tilde\Ps^n(0)=\tilde{X}^n(0)\in\R^I_+$ and,
from~(\ref{417}), that
$(\Ps^n, \tilde\Ps^n)\in\calD(2(2+\kappa_4))$ for all large $n$.
Since $0\leq B^n_i(s)\leq1$, $T^n_i(s)\in[0, \tau_n]$ for all $s\in
[0,\tau_n]$.
Hence from (\ref{411}) we have for $(\tilde{A}^n, \tilde{S}^n)\in\A_r$
\[
\sup_{[0,\tau_n]}\bigl|\tilde\psi^2_i(
\rho_it)-\tilde {S}^n_i\bigl(T^n_i(t)
\bigr)\bigr|\leq r_1+\sup_{[0,\tau_n]}\bigl|\tilde
\psi^2_i(\rho _it)-\tilde{
\psi}^2\bigl(T^n_i(t)\bigr)\bigr|.
\]
Again using the continuity of $\tilde\psi^2$,
we can choose $r_2>0$ small enough such that
$\osc_{r_2}[\tilde\psi^2]<r_1$.
Since $\frac{b_n}{\sqrt{n}}\to0$, we note from (\ref{417}) that for
all large $n$, and
all $i$, $\sup_{[0,\tau_n]}|\rho_i t-T^n_i(t)|<r_2$.
Since $\tilde X^n(0)\to x$, $y^n\to y$ and $\theta^n\to\theta$, it
follows that
\[
\bigl|\tilde\Ps^n_i-\tilde\psi^2_i(
\rho_i\cdot)\bigr|^*_{\tau_n}<3r_1
\]
for all large $n$.
Now taking $\kappa=2(2+\kappa_4)$, we choose $r_1$ sufficiently small
[see~(\ref{315})] so that for all $n$ large we have
\[
\bigl\|\bzeta[\tilde\psi]-\hat{\bzeta}\bigl[\Ps^n,\tilde\Ps^n
\bigr]\bigr\|^*_{\tau
_n} \leq\eps.
\]
Now choosing $r<\eps/(3\sqrt{I})$ and using \eqref{12} and \eqref{13},
for $(\tilde{A}^n, \tilde{S}^n)\in\A_r$ and all large $n$, we have
%
\begin{equation}
\bigl\|\varphi-\Ph^n\bigr\|^*_{\tau_n}\leq4\eps.\label{419}
\end{equation}
Let $\kappa_5=(\|\varphi\|^*+4)$.
Denote by $\omega_h$ [resp., $\omega_g$] the modulus of continuity of
$h$ [resp., $g$]
over $\{q\dvtx \|q\|\leq\kappa_5 \}$.
Then by \eqref{14}, on the indicated event, for all large $n$,
\[
\int_0^{\tau_n}h\bigl(\tilde{X}^n(s)
\bigr)\,ds \geq\int_0^{\tau_n}h\bigl(
\Ph^n(s)\bigr)\,ds \geq\int_0^{\tau_n} h
\bigl(\varphi(s)\bigr)\,ds-T\omega_h(4\eps).
\]
Combined with \eqref{444} this gives
\[
\int_0^Th\bigl(\tilde X^n(s)
\bigr)\,ds\ge\int_0^Th\bigl(\varphi(s)\bigr)\,ds-T
\omega _h(4\eps).
\]
A similar argument gives
\[
g\bigl(\tilde X^n(T)\bigr)= g\bigl(\varphi(T)\bigr)
\chi_{\{T\leq\tau_n\}}+g\bigl(\varphi(T)\bigr)\chi_{\{T>\tau
_n\}} \ge g\bigl(
\varphi(T)\bigr)-\omega_g(4\eps).
\]
Hence for all large $n$,
\begin{eqnarray*}
\E\bigl[e^{b_n^2[\int_0^Th(\tilde{X}^n(s))\,ds+g(\tilde{X}^n(T))]}\bigr] &\geq& \E \bigl[e^{b_n^2[\int_0^Th(\tilde{X}^n(s))\,ds+g(\tilde{X}^n(T))]}
\chi_{\{(\tilde{A}^n, \tilde{S}^n)\in\A_r\}} \bigr]
\\
&\ge& \E \bigl[e^{b_n^2[\int_0^Th(\varphi(s))\,ds+g(\varphi(T))
-a(\eps)]}\chi_{\{(\tilde{A}^n, \tilde{S}^n)\in\A_r\}} \bigr],
\end{eqnarray*}
where $a(\eps)=[T\omega_h(4\eps)+\omega_g(4\eps)]\to0$ as $\eps
\to0$.
We now use Assumption \ref{moderate}. Since $\A_r$ is open,
\[
\p\bigl(\bigl(\tilde{A}^n, \tilde{S}^n\bigr)\in
\A_r\bigr)\geq e^{-b_n^2[\inf_{\psi
\in\A_r}\Ir(\psi)+\eps]}\geq e^{-b_n^2[\Ir(\tilde\psi)+\eps]}
\]
holds for all sufficiently large $n$.
Hence we have from \eqref{75} and (\ref{412}) that for all large $n$,
\begin{eqnarray*}
V^n\bigl(\tilde{X}^n(0)\bigr)+\eps&\geq& J\bigl(
\tilde{X}^n(0), B^n\bigr)
\\
&\geq& \int_0^Th\bigl(\varphi(s)\bigr)\,ds+g
\bigl(\varphi(T)\bigr)-\Ir(\tilde\psi )-a(\eps)-\eps.
\end{eqnarray*}
Therefore
\[
\liminf_{n\to\iy} V^n\bigl(\tilde{X}^n(0)
\bigr)\geq c\bigl(\tilde\psi, \bzeta [\tilde\psi]\bigr)-a(\eps)-2\eps,
\]
and letting $\eps\to0$, we obtain
$
\liminf_{n\to\iy} V^n(\tilde{X}^n(0))\geq c(\tilde\psi, \bzeta
[\tilde\psi])$.
Finally, since $\tilde\psi\in\calP$ is arbitrary we have from (\ref
{305}) that
$
\liminf_{n\to\iy} V^n(\tilde{X}^n(0))\geq V(x)$.
\end{pf*}

\subsection{Upper bound}\label{sec42}

\begin{theorem}\label{th42}
Let Assumptions \ref{moderate} and \ref{mini} hold. If either $g$ or
$h$ is unbounded,
let also Assumption \ref{unbounded} hold.
Then $\limsup V^n(\tilde{X}^n(0))\leq V(x)$.
\end{theorem}

The proof is based on the construction and analysis of a particular policy,
described below in equations \eqref{20}--\eqref{25}.
To see the main idea behind the structure of the policy, refer to
equations \eqref{eqn2} and \eqref{06},
which describe the dependence of
the scaled process $\tilde X^n$ on the stochastic primitives $\tilde
A^n$, $\tilde S^n$
and the control process $B^n$ [recall from \eqref{02} that $T^n$ is an
integral form of~$B^n$].
Because of the amplifying factor\vadjust{\goodbreak} $\sqrt n/b_n$ which appears in the
expression \eqref{06} in front of
\[
\rho_it-T^n_i(t)=\int_0^t
\bigl(\rho_i-B^n_i(s)\bigr)\,ds,
\]
it is seen that fluctuations of $B^n$ about its center $\rho$,
at scale as small as $b_n/\sqrt n$, cause order-one displacements in
$\tilde X^n$.
Initially, the policy drives the process $\tilde X^n$ from
the initial position $\tilde X^n(0)\approx x$ to the corresponding
point on the minimizing curve, $f(\theta\cdot x)$, in a short time.
This is reflected in the choice of the constant $\ell$ applied during
the first time interval $[0,v)$; see the first line of \eqref{21}.
Afterwards, the policy mimics the behavior of the optimal strategy for
the game, namely $\hat{\bzeta}$.
This is performed by applying $F^n$; see the third line of \eqref{21},
which consists of the response of $\hat{\bzeta}$, in differential
form, to the stochastic data $P^n$; see \eqref{23}.

\begin{pf*}{Proof of Theorem \ref{th42}}
Given a constant $\Del$, define
%
\begin{equation}
\label{60} \calD_\Del=\bigl\{\psi\in\calD\bigl([0, T],
\R^{2I}\bigr) \dvtx \Ir(\psi)\leq \Del\bigr\}.
\end{equation}
By the definition of $\Ir$ (from Section~\ref{sec2}), $\calD_\Del$
is a compact set containing absolutely
continuous paths starting from zero (particularly, $\calD_\Del\subset
\calP$),
with derivative having $L^2$-norm uniformly bounded.
Consequently, for a constant $M=M_\Del$ and all $\psi\in\calD_\Del
$, one has
$\|\psi^1\|^*+\|\psi^2\|^*\leq M$. Consider the set $\calD(M+1)$
[see \eqref{19}],
let $\eps\in(0,1)$ be given, and choose $\delta_1, \delta, \eta>0$
as in (\ref{315}) and (\ref{316}), corresponding to $\eps$ and
$\kappa=M+1$.
Assume, without loss of generality, that $\delta_1\vee\delta<\eps$.
It follows from the $L^2$ bound alluded to above, that for each fixed
$\Del$,
the members of $\calD_\Del$ are equicontinuous. Hence one can choose
$v_0\in(0,\eta)$ (depending on
$\Del$), such that
%
\begin{equation}
\osc_{v_0}\bigl(\psi^l_i\bigr)<
\frac{\delta_1\wedge\delta}{4\sqrt{2I}}\qquad \mbox{for all } \psi=\bigl(\psi^1,
\psi^2\bigr)\in\calD_\Del, l=1,2, i\in\I.\label{456}
\end{equation}
Recall from \eqref{74}--\eqref{76} the notation $d$, $\Ups$ and $\|
\cdot\|^\circ$.
As in the proof of Theorem \ref{th41}, we set for $\tilde\psi\in
\calP$,
\[
\A_r(\tilde\psi)=\bigl\{\psi\in\calD\bigl([0,T], \R^{2I}
\bigr) \dvtx d(\psi,\tilde\psi)<r\bigr\}.
\]
Noting that, for any $f\in\Ups$,
\begin{eqnarray*}
\bigl\|\psi(t)-\tilde\psi(t)\bigr\| &\leq& \bigl\|\psi(t)-\tilde\psi\bigl(f(t)\bigr)\bigr\|+\bigl\|
\tilde\psi\bigl(f(t)\bigr)-\tilde\psi(t)\bigr\|,
\\
\bigl|f(\cdot)-\cdot\bigr|^*_T &\leq& T\bigl(e^{\|f\|^\circ}-1\bigr),
\end{eqnarray*}
it follows by equicontinuity that it is possible to choose
$v_1>0$ such that, for any $\tilde\psi\in\calD_\Del$,
%
\begin{equation}
\psi\in\A_{v_1}(\tilde\psi) \mbox{ implies } \|\psi-\tilde\psi\|^*<
\frac{\delta_1}{4}.\label{457}
\end{equation}
Let $v_2=\min\{v_0, v_1, \frac{\eps}{2}\}$. Since $\calD_\Del$ is
compact and $\Ir$ is lower semicontinuous,
one can find a finite number of members $\psi^1, \psi^2,\ldots,$
$\psi^N$
of $\calD_\Del$, and positive constants $v^1,\ldots, v^N$ with
$v^k<v_2$, satisfying
$\calD_\Del\subset\bigcup_k \A^k$, and
%
\begin{equation}
\inf\bigl\{\Ir(\psi)\dvtx \psi\in\oo{\A^k}\bigr\}\geq\Ir\bigl(
\psi^k\bigr)-\frac{\eps}{2},\qquad k=1,2,\ldots,N, \label{421}
\end{equation}
where, throughout, $\A^k:=\A_{v^k}(\psi^k)$.

We next define a policy for which we shall prove that the lower bound is
asymptotically attained.
Fix $n\in\N$.
Recall \eqref{01}, \eqref{02} and \eqref{03} by which
%
\begin{equation}
\label{20} \cases{ D^n_i=S^n_i
\circ T^n_i,\vspace*{2pt}
\cr
\displaystyle T^n_i=
\int_0^\cdot B^n_i(s)\,ds,
\vspace*{2pt}
\cr
X^n_i=X^n_i(0)+A^n_i-D^n_i.}
\end{equation}
Recall the scaled processes \eqref{22} and let also
%
\begin{equation}
\label{23} \cases{ \tilde D^n_i=\tilde
S^n_i\circ T^n_i,\vspace*{2pt}
\cr
P^n=\bigl(\tilde A^n,\tilde D^n\bigr). }
\end{equation}
The analogy between the queueing system dynamics \eqref{eqn2} and the
game dynamics~\eqref{201}
suggests that the policy should be designed in such a way that
$\frac{\mu_i\sqrt n}{b_n}\int_0^\cdot(\rho_i-B^n_i(s))\,ds\approx
\bzeta_i[P^n]$
holds for each $i$.
Equivalently, one should have
$\int_0^t B^n_i(s)\,ds\approx\rho_it-\frac{b_n}{\mu_i\sqrt n}\bzeta
_i[P^n](t)$.
A straightforward discretization approach fails to provide an
admissible control.
A version of this approximate equality that
does define an admissible control is as follows.
Denote
%
\begin{equation}
\label{79} \Th(a, b)=a\chi_{[0,1]}(a)\chi_{[0,1]}(b), \qquad a, b\in
\R.
\end{equation}
Let $\ell=f(x\cdot\theta)-x$ and $v=\frac{v_2}{2}\wedge\frac{T}{4}$.
For $i\in\I$, assume that $B^n_i$ is given by
%
\begin{equation}
\label{24} B^n_i(t)=C^n_i(t)
\chi_{\{X^n_i(t)>0\}},\qquad t\in[0,T],
\end{equation}
where, for $t\in[0,T]$,
%
\begin{equation}\qquad
\label{21} C^n_i(t)=\cases{ \displaystyle\Th \Biggl(
\rho_i-\frac{b_n}{\mu_i\sqrt n}\frac{\ell_i}{v}, \sum
_{k=1}^I \biggl(\rho_k-
\frac{b_n}{\mu_k\sqrt n}\frac{\ell_k}{v} \biggr)^+ \Biggr),\vspace*{2pt}\cr
\qquad\hspace*{13.5pt} \mbox{if } t\in[0, v),
\vspace*{2pt}
\cr
\rho_i, \qquad  \mbox{if } t\in[v, 2v),\vspace*{2pt}
\cr
\displaystyle\Th \Biggl(\rho_i-F^n_i(t-v), \sum
_{k=1}^I\bigl(\rho_k-F^n_k(t-v)
\bigr)^+ \Biggr),\vspace*{2pt}\cr
\hspace*{37pt}\mbox{if } \bigl\|P^n\bigr\|^*_{t-v}< M+2, t
\in\bigl[jv,(j+1)v\bigr), j=2,3,\ldots,\vspace *{2pt}
\cr
\rho_i,\qquad \mbox{if }
\bigl\|P^n\bigr\|^*_{t-v}\ge M+2, t\in\bigl[jv,(j+1)v\bigr), j=2,3,\ldots,}
\end{equation}
and
%
\begin{eqnarray}
\label{25} F^n_i(u)=\frac{b_n}{\mu_i\sqrt{n}}
\frac{\hat{\bzeta}_i[P^n](jv)-
\hat{\bzeta}_i[P^n]((j-1)v)}{v},
\nonumber
\\[-8pt]
\\[-8pt]
 \eqntext{u\in\bigl[jv,(j+1)v\bigr), j=1,2,\ldots.}
\end{eqnarray}
Let us argue that these equations uniquely define a policy.
To this end, consider equations \eqref{20}, \eqref{23}, \eqref{24},
\eqref{21}, \eqref{25},
along with the obvious relations between scaled and unscaled processes,
as a set of equations for $X^n,D^n,T^n,P^n,B^n,C^n,F^n$ (and the scaled
versions $\tilde X^n,\tilde D^n$),
driven by the data $(A^n,S^n)$ [equivalently, $(\tilde A^n,\tilde
S^n)$], and satisfying
the initial condition $X^n(0)$. Arguing by induction on the jump times
of the processes
$A^n$ and $S^n$, and using the causality of the map $\hat{\bzeta}$,
it is easy
to see that this
set of equations has a unique solution.
Moreover, this solution is consistent with the model equations
\eqref{01}--\eqref{03}.
The processes alluded to above are therefore well defined.

We now show that $B^n\in\Bf^n$.
To see that $B^n$ has RCLL sample paths, note first that, by construction,
$F^n$, $X^n$ are piecewise constant with finitely many jumps, locally,
hence so is $B^n$.
Therefore the existence of left limits follows. Right continuity
follows from the fact
that $X^n$, $F^n$ and consequently $C^n$ have this property.
The other elements in the definition of an admissible control hold by
construction.
Thus $B^n\in\Bf^n$ for $n\in\N$. As a result,
%
\begin{equation}
\label{45} V^n\bigl(\tilde X^n(0)\bigr)\le
J^n\bigl(\tilde X^n(0),B^n\bigr).
\end{equation}

Our convention in this proof will be that $c_1,c_2,\ldots$ denote
positive constants that do not depend on
$n,\eps,v,\Del$. Also, the notation \eqref{70}--\eqref{83} will be
used extensively.

Let, for $k=1,\ldots,N$,
\[
\bigl(\varphi^k,\xi^k,\zeta^k,
w^k\bigr)=\bigl(\bph\bigl[\psi^k\bigr],\bxi\bigl[
\psi^k\bigr],\bzeta \bigl[\psi^k\bigr],\bw\bigl[
\psi^k\bigr]\bigr).
\]
Write $\psi^k$ as $(\psi^{k,1},\psi^{k,2})$.
Note that $\varphi^k$ is the dynamics corresponding to
$\psi^k$ and~$\zeta^k$.
Let $\La^n=\|\tilde{A}^n\|^*_T+\|\tilde{S}^n\|^*_T$, and define
%
\begin{equation}
\label{2222} \Om^n_k=\bigl\{\bigl(\tilde{A}^n,
\tilde{S}^n\bigr)\in\A^k\bigr\}, \qquad k=1,\ldots, N.
\end{equation}
We prove the result in number of steps.
In steps 1--4 we shall show that for a constant $c_1$, for all $n\ge
n_0(\eps,v)$,
%
\begin{equation}
\label{41} \bigl\|\tilde X^n\bigr\|^*_T\le c_1
\bigl(1+\La^n\bigr)
\end{equation}
and
%
\begin{equation}
\label{42} \sup_{[v,T]}\bigl\|\tilde X^n-
\varphi^k\bigr\|\le c_1\eps \qquad\mbox{on } \Om
^n_k, k=1,2,\ldots,N.
\end{equation}
The final step will then use these estimates to conclude the result.

\textit{Step} 1: The goal of this step is to show \eqref{425}
below which is the key estimate in proving \eqref{41}.
By Proposition \ref{prop2}(ii),
%
\begin{equation}
\label{43} \bigl\|\hat{\bzeta}\bigl[P^n\bigr]\bigr\|^*_t\le
c_2\bigl(1+\bigl\|P^n\bigr\|^*_t\bigr).
\end{equation}
Therefore
%
\begin{equation}
\bigl\|F^n\bigr\|^*_t\leq \frac{b_n}{\sqrt{n}}\frac{c_3}{v}
\bigl(1+\bigl\|P^n\bigr\|^*_t\bigr). \label{423}
\end{equation}
Since $\rho_i\in(0,1)$ for all $i\in\I$, we note from (\ref{423})
that for all sufficiently large $n$, for any $t\in[2v, T]$,
\[
\bigl\|P^n\bigr\|^*_{t-v}< M+2 \mbox{ implies } \sum
_i\bigl(\rho_i-F^n_i(t-v)
\bigr)^+=\sum_i\bigl(\rho_i-F^n_i(t-v)
\bigr)\leq1
\]
as $\sum_{i}F^n_i(u)\geq0$ for all $u\in[v, T]$.
Define
\[
\hat\tau_n=\inf\bigl\{t\geq0 \dvtx\bigl \|P^n(t)\bigr\|\geq M+2
\bigr\}.
\]
It is easy to check by definition of $C^n_i$, and using the fact $\rho
_i\in(0,1)$ and the
convergence $b_n/\sqrt{n}\to0$,
that for all large $n$, on the event $\{\hat\tau_n\leq v\}$,
\[
\sup_{t\in[0,T]}\frac{\sqrt{n}}{b_n} \biggl|\rho_i t-\int
_0^t C^n_i(s)\,ds \biggr|\leq
c_4.
\]
Next consider the event $\{\hat\tau_n>v\}$. Using \eqref{79}, \eqref
{21} and
\eqref{423}, one has for all sufficiently large $n$,
%
\begin{equation}
\label{48} C^n_i(t)=\cases{
\displaystyle \rho_i-\frac{b_n}{\mu_i\sqrt n}\frac{\ell_i}{v}, &\quad
$\mbox{if } t\in[0,v),$
\vspace*{2pt}\cr
\rho_i, & \quad$\mbox{if } t\in[v,2v),$
\vspace*{2pt}\cr
\rho_i-F^n_i(t-v), & \quad$\mbox{if } t\in[2v,
\hat{\tau}_n+v),$
\vspace*{2pt}\cr
\rho_i, & \quad$\mbox{if } t\in[\hat{\tau}_n+v,T].$}
\end{equation}
Thus, on $\{\hat\tau_n>v\}$,
\[
\sup_{t\in[0,2v]} \biggl|\rho_i t-\int_0^t
C^n_i(s)\,ds \biggr|\leq c_5\frac{b_n}{\sqrt n},
\]
while
%
\begin{equation}
\label{44} \sup_{t\in[2v,T]} \biggl|\rho_i t-\int
_0^t C^n_i(s)\,ds \biggr|\leq
c_5\frac{b_n}{\sqrt n}+\sup_{t\in[2v,\hat\tau_n+v]} \biggl|\int
_{2v}^t F^n_i(s-v)\,ds \biggr|.
\end{equation}
Consider $j\geq2$ and $jv\leq t< (j+1)v$. Then by the definition of $F^n$,
%
\begin{eqnarray}
\label{46}
\nonumber
\int_{2v}^tF^n_i(s-v)\,ds
&=& \int_{2v}^{jv}F^n_i(s-v)\,ds
+ \int_{jv}^tF^n_i(s-v)\,ds
\\
&= &\frac{b_n}{\mu_i\sqrt{n}} \bigl[\hat{\bzeta}_i
\bigl[P^n\bigr]\bigl((j-2)v\bigr)-\hat{\bzeta}_i
\bigl[P^n\bigr](0)\bigr]
\\
&&{} + \frac{b_n}{\mu_i\sqrt{n}}\frac{t-jv}{v} \bigl[\hat{\bzeta}_i
\bigl[P^n\bigr]\bigl((j-1)v\bigr)-\hat{\bzeta}_i
\bigl[P^n\bigr]\bigl((j-2)v\bigr)\bigr].\nonumber
\end{eqnarray}
Combining this identity with \eqref{43}
shows that the last term on \eqref{44} is bounded by
\[
\sup_{t\in[2v,\hat\tau_n+v]}\frac{b_n}{\mu_i\sqrt n}4c_2\bigl(1+
\bigl\|P^n\bigr\| ^*_{t-v}\bigr)\le \frac{b_n}{\mu_i\sqrt n}4c_2
\bigl(1+\La^n\bigr),
\]
where in the last inequality we also used the fact that $T^n_i(t)\le
t$, by which
$|\tilde D^n_i|^*_t=|\tilde S^n_i\circ T^n_i|^*_t\le|\tilde S^n_i|^*_t$.
We conclude that, for all sufficiently large $n$,
%
\begin{equation}
\sup_{t\in[0,T]}\frac{\sqrt{n}}{b_n}\biggl |\rho_i t-\int
_0^t C^n_i(s)\,ds \biggr|\leq
c_6\bigl(1+\La^n\bigr). \label{425}
\end{equation}

\textit{Step} 2: We prove \eqref{41}.
The argument is based on the \textit{Skorohod problem} (see, e.g., \cite{chen-yao}) and the estimate \eqref{425}.
To this end, rewrite (\ref{eqn2}) as $\tilde X^n_i=\hat Y^n_i+\hat
Z^n_i$, where
\begin{eqnarray*}
\hat Y^n_i(t) &=& \tilde{X}^n_i(0)+y^n_it+
\tilde{A}^n_i(t)-\tilde {S}^n_i
\bigl(T^n_i(t)\bigr) +\frac{\mu^n_i}{n}
\frac{\sqrt{n}}{b_n} \biggl(\rho_i t-\int_0^tC^n_i(s)\,ds
\biggr),
\\
\hat Z^n_i(t) &=& \frac{\mu^n_i}{n}\frac{\sqrt{n}}{b_n}
\int_0^tC^n_i(s)
\chi_{\{\tilde{X}^n_i(s)=0\}}\,ds.
\end{eqnarray*}
Since for each $i$, $\tilde X^n_i$ is nonnegative and $\hat Z^n_i$ is
nonnegative, nondecreasing
and increases only
when $\tilde X^n_i$ is equal to zero, it follows that $(\tilde
X^n_i,\hat Z^n_i)$ is the solution to the
Skorohod problem for data $\hat Y^n_i$; see \cite{chen-yao} and
\cite{cheman} for this well-known characterization of the Skorohod map
\eqref{skoro}.
As a result, for all large $n$,
%
\begin{equation}
\label{26} \bigl|\hat Z^n_i\bigr|^*_T+\bigl|\tilde
X^n_i\bigr|^*_T\leq4\bigl|\hat Y^n_i\bigr|^*_T
\le c_7\bigl(1+\La^n\bigr),
\end{equation}
where we used (\ref{425}) and the convergence of $\mu^n_i/n$, $\tilde
X^n_i(0)$ and $y^n_i$.
This shows~\eqref{41}.

\textit{Step} 3: Here we analyze the events $\Om^n_k$, showing
that on these events one has,
for large $n$, that $\mu_i\frac{\sqrt n}{b_n}(\rho_it-\int_0^tC^n_i(s)\,ds)$ is close to $\zeta^k_i$.
First, using
\[
\rho_it-T^n_i(t)=\rho_it-\int
_0^tC^n_i(s)\,ds+\int
_0^tC^n_i(s)\chi
_{\{\tilde C^n_i(s)=0\}}\,ds,
\]
we obtain from \eqref{425} and \eqref{26}, for all large $n$,
\[
\sup_{t\in[0, T]}\frac{\mu^n_i}{n}\frac{\sqrt{n}}{b_n}\bigl|
\rho_it-T^n_i(t)\bigr| \leq c_8
\bigl(1+\La^n\bigr).
\]
Therefore we obtain that, for all large $n$, on the event $\bigcup_k\Om^n_k$,
%
\begin{equation}
\sup_{t\in[0, T]}\bigl|\rho_it-T^n_i(t)\bigr|
\leq\frac{v}{2}. \label{427}
\end{equation}
This shows that under the policy $B^n$, on $\bigcup_k\Om^n_k$,
the average effort given by the server to class-$i$ customers
is equal to $\rho_i$ asymptotically.
Abusing the notation and writing
$\psi^{k,2}(T^n(\cdot))$ for
$(\psi^{k,2}_1(T^n_1(\cdot)),\ldots,\psi^{k,2}_{I}(T^n_I(\cdot)))$,
using (\ref{456}) and (\ref{427}) for the choice of $v$, we have
%
\begin{equation}
\label{2202} \sup_{t\in[v, T]}\bigl\|\psi^{k,2}
\bigl(T^n(t)\bigr)-\brho\bigl[\psi^{k,2}\bigr](t-v)\bigr\|\leq
\Biggl[\sum_{i=1}^I \bigl(
\operatorname{osc}_{2v}\bigl(\psi^{k, 2}_i\bigr)
\bigr)^2 \Biggr]^{1/2}\leq\frac{\delta_1}{4},
\end{equation}
on $\Om^n_k$, for all $n$ large.

Next, we estimate
$\tilde{S}^n(T^n(t))-\brho[\psi^{k,2}](t-v)$ on $\Om^n_k$.
Using \eqref{457}, for all large~$n$,
%
\begin{eqnarray}
\label{428}
\nonumber
&&\sup_{t\in[v, T]}\bigl\|\tilde{S}^n
\bigl(T^n(t)\bigr) -\brho\bigl[\psi^{k,2}\bigr](t-v)\bigr\|
\\
&&\qquad\le\bigl\|\tilde{S}^n\bigl( T^n(\cdot)\bigr)-
\psi^{k,2}\bigl(T^n(\cdot)\bigr)\bigr\|^*+\sup
_{t\in[v, T]}\bigl\|\psi^{k,2}\bigl(T^n(t)\bigr) -
\brho\bigl[\psi^{k,2}\bigr](t-v)\bigr\|
\\
&&\qquad\le\frac{\delta_1}{4}+\frac{\delta_1}{4}=\frac{\delta_1}{2},\nonumber
\end{eqnarray}
where for the first estimate we have used (\ref{457}) and for second
we have used~\eqref{2202}.

Finally, we show the two estimates \eqref{47} and \eqref{430}, below.
Note that on $\Om^n_k$ one has $\hat\tau_n\ge T$ for all large $n$
[as follows
by $\|P^n\|^*_T=\|\tilde A^n\|^*_T+\|\tilde D^n\|^*_T\le
\|\tilde A^n\|+\|\tilde S^n\|<M+2$ by the discussion in the beginning
of the proof \eqref{457}].
As a result, \eqref{48} is applicable. In particular, for all large $n$,
%
\begin{equation}
\label{47} \mu_i\frac{\sqrt{n}}{b_n} \biggl(\rho_it-
\int_0^tC^n_i(s)\,ds
\biggr)-\frac{t}{v}\ell_i=0,\qquad  t\in[0,v).
\end{equation}

Now for $k=1,2,\ldots,N$, consider
\[
\hat W_{i,k}^n(t):=\mu_i\frac{\sqrt{n}}{b_n}
\biggl(\rho_it-\int_0^tC^n_i(s)\,ds
\biggr) -\zeta_i^k(t-v),\qquad  t\in[v,T],
\]
on the event $\Om^n_k$. We note from (\ref{73}) that $\zeta
^k(0)=\ell$. Hence
for $t\in[v,2v)$ and all large $n$, we have from (\ref{456}) and
(\ref{316}) that
\[
\bigl|\hat W_{i,k}^n(t)\bigr| =\bigl|\ell_i-
\zeta^k_i(t-v)\bigr|\leq\eps.
\]
Next consider $t\in[2v, T]$ and integer $j$ for which $jv\leq
t<(j+1)v$. From calculation \eqref{46},
for large $n$,
\begin{eqnarray*}
\mu_i\frac{\sqrt{n}}{b_n} \biggl(\rho_it-\int
_0^t C^n_i(s)\,ds \biggr)
&=& \ell_i+ \mu_i\frac{\sqrt{n}}{b_n}\int
_{2v}^t F^n_i(s-v)\,ds
\\
&=&\hat{\bzeta}_i\bigl[P^n\bigr]\bigl((j-2)v\bigr)
\\
&&{} +\frac{t-jv}{v}\bigl[\hat{\bzeta}_i\bigl[P^n
\bigr]\bigl((j-1)v\bigr)-\hat{\bzeta }_i\bigl[P^n\bigr]
\bigl((j-2)v\bigr)\bigr].
\end{eqnarray*}
Hence
\begin{eqnarray*}
\bigl|\hat W_{i,k}^n(t)\bigr|&\leq&\bigl|\hat{\bzeta}_i
\bigl[P^n\bigr]\bigl((j-2)v\bigr)-\zeta^k_i(t-v)\bigr|\\
&&{}+
\bigl|\hat{\bzeta}_i\bigl[P^n\bigr]\bigl((j-1)v\bigr) -\hat{
\bzeta}_i\bigl[P^n\bigr]\bigl((j-2)v\bigr)\bigr|.
\end{eqnarray*}
For large $n$,
\begin{eqnarray*}
&&\bigl|\hat{\bzeta}_i\bigl[P^n\bigr]\bigl((j-2)v\bigr) -
\zeta^k_i(t-v)\bigr|
\\
&&\qquad\leq\bigl|\hat{\bzeta}_i\bigl[P^n\bigr]\bigl((j-2)v\bigr) -
\hat{\bzeta}_i\bigl[\psi^{k,1},\psi^{k,2}\circ
T^n\bigr]\bigl((j-2)v\bigr)\bigr|
\\
& &\qquad\quad{}+\bigl|\hat{\bzeta}_i\bigl[\psi^{k,1},\psi^{k,2}
\circ T^n\bigr]\bigl((j-2)v\bigr) -\hat{\bzeta}_i\bigl[
\psi^{k,1},\brho\bigl[\psi^{k,2}\bigr]\bigr]\bigl((j-2)v\bigr)\bigr|
\\
&&\qquad\quad{} +\bigl|\zeta^k_i\bigl((j-2)v\bigr)-\zeta^k_i(t-v)\bigr|
\\
&&\qquad \leq3\eps,
\end{eqnarray*}
where the first quantity is estimated using (\ref{457}) and~(\ref
{315}), the second using
(\ref{427}) and~(\ref{315}), and the third using (\ref{456}) and
(\ref{316}).
A similar estimate gives, for all large~$n$,
\[
\bigl|\hat{\bzeta}_i\bigl[P^n\bigr]\bigl((j-1)v\bigr)-\hat{
\bzeta}_i\bigl[P^n\bigr]\bigl((j-2)v\bigr)\bigr|\le3\eps.
\]
Hence for all large $n$, on $\Om^n_k$,
%
\begin{equation}
\sup_{t\in[v,T]}\bigl|\hat W_{i,k}^n(t)\bigr|\le6\eps.
\label{429}
\end{equation}
Using (\ref{429}) and (\ref{425}), for all large $n$, on $\Om^n_k$,
%
\begin{equation}
\sup_{t\in[v,T]} \biggl|\frac{\mu^n_i}{n}\frac{\sqrt{n}}{b_n} \biggl(
\rho_it-\int_0^tC^n_i(s)\,ds
\biggr) -\zeta^k_i(t-v) \biggr|\le7\eps. \label{430}
\end{equation}
Thus we see from \eqref{eqn2}, \eqref{201}, \eqref{428} that under
the defined policy $B^n$ the scaled process $\tilde X^n$ stays near the
path $\varphi^k$ on $\Om^n_k$ provided we can
control the error that arises from the server idleness. In the next
step we show that this can be done.

\textit{Step} 4: Now we prove \eqref{42}.
Recall $\varphi^k=\bph[\psi^k]$.
The goal of this step is to estimate the difference between $\tilde
X^n$ and $\varphi^k$ on $\Om^n_k$.
To this end, let first
\[
\tilde\varphi^k(t)=\cases{ %
\displaystyle x+
\frac{t}{v}\ell,& \quad $\mbox{for } t\in[0,v),$
\vspace*{2pt}\cr
\varphi^k(t-v), & \quad $\mbox{for } t\in[v,T].$}
\]
Recall from step 2 that $\tilde X^n_i$ solves the Skorohod problem for
$\hat Y^n_i$. Note also that $\tilde\varphi^k_i\ge0$.
Thus using the Lipschitz property of the Skorohod map we have on $\Om^n_k$
%
\begin{equation}
\label{50} \bigl|\tilde{X}^n_i-\tilde
\varphi^k_i\bigr|^*_T\leq2\bigl|\hat
Y^n_i-\tilde \varphi^k_i\bigr|^*_T.
\end{equation}
For $t\in[0, v]$ and $n$ large, we have, using the definition of $\hat
Y^n$ and \eqref{47},
%
\begin{eqnarray}
\label{51}
&&\bigl|\hat Y^n_i(t)-\tilde
\varphi^k_i(t)\bigr|
\nonumber\\
&&\qquad\leq\bigl|\tilde{X}^n_i(0)-x_i\bigr|+v\bigl|y^n_i\bigr|
+ \bigl|\tilde{A}^n_i(t) -\tilde{S}^n_i
\bigl(T^n_i(t)\bigr)\bigr|
\nonumber
\\[-8pt]
\\[-8pt]
\nonumber
&&\qquad\quad{}+ \biggl|\frac{\mu^n_i}{n}-
\mu_i \biggr| \biggl|\frac{\sqrt{n}}{b_n} \biggl(\rho_i t-\int
_0^tC^n_i(s)\,ds \biggr)
\biggr|
\\
&&\qquad\leq c_9\eps\nonumber
\end{eqnarray}
on $\Om^n_k$, where we use (\ref{456}), (\ref{457}) and (\ref{425}).
Moreover, for $t\in[v,T]$, by the definition of $\hat Y^n$ and $\tilde
\varphi^k$,
\begin{eqnarray*}
\hat Y^n_i(t)-\tilde\varphi^k_i(t)
&=& \tilde X^n_i(0)+y^n_it+
\tilde A^n_i(t)-\tilde S^n_i
\bigl(T^n_i(t)\bigr)\\
&&{}+\frac
{\mu^n_i}{n}\frac{\sqrt n}{b_n}
\biggl(\rho_it-\int_0^tC^n_i(s)\,ds
\biggr)
\\
&&{} -\zeta^k_i(t-v)-x_i-y_i(t-v)-
\psi^{k,1}_i(t-v)+\brho_i\bigl[
\psi^{k,2}\bigr](t-v).
\end{eqnarray*}
Hence, using (\ref{456}), (\ref{457}), (\ref{428}) and (\ref{430}),
estimate \eqref{51}
is valid for $t\in[v, T]$ as well. Namely, $|\hat Y^n_i-\tilde\varphi
^k_i|^*_T\le c_9\eps$ on $\Om^n_k$
for large $n$. Thus using \eqref{50},
$\|\tilde X^n-\tilde\varphi^k\|^*\le c_{10}\eps$ on $\Om^n_k$ for
large $n$.
By the definition of $\tilde\varphi^k$ and \eqref{73}, \eqref{316},
\eqref{456} we obtain that, for all sufficiently large $n$,
\eqref{42} holds.

\textit{Step} 5:
Finally, in this step, we rely on property \eqref{421} to complete the proof.
Since $\varphi^k$ is bounded, and so is $\tilde X^n$ on $\Om^n_k$, it
follows from \eqref{42}
by continuity of $h$ and $g$ that, for all large $n$, on $\Om^n_k$,
%
\begin{equation}
\label{49} \biggl|\int_0^Th\bigl(
\varphi^k(s)\bigr)\,ds+g\bigl(\varphi^k(T)
\bigr)-H^n \biggr|\le\omega (\eps),
\end{equation}
where
\[
H^n=\int_0^Th\bigl(
\tilde{X}^n(s)\bigr)\,ds+g\bigl(\tilde{X}(T)\bigr),
\]
and $\omega=\omega_\Del$ satisfies $\omega(a)\to0$ as $a\to0$,
for any $\Del$.
By \eqref{41} and the growth condition on $h$ and $g$, $H^n\le
c_{11}(1+\La^n)$.
Hence given any $\Del_1>0$,
\[
H^n> \Del_1 \mbox{ implies } \La^n >
c_{11}^{-1}\Del_1-1=:G(\Del_1).
\]
Therefore
\begin{eqnarray}
\label{432} \E\bigl[e^{b_n^2H^n}\bigr] &\leq& \E\bigl[e^{b_n^2[H^n\wedge\Del_1]}
\bigr]+
\E\bigl[e^{b_n^2H^n} \chi_{\{H^n>\Del_1\}}\bigr]
\nonumber
\\[-8pt]
\\[-8pt]
\nonumber
&\leq& \E\bigl[e^{b_n^2[H^n\wedge\Del_1]}\bigr] + \E\bigl[e^{b_n^2c_{11}(1+\La^n)}
\chi_{\{\La^n>G(\Del_1)\}}\bigr].
\end{eqnarray}
Now we estimate both terms on the RHS of (\ref{432}).
Denote $\B=(\bigcup_{k=1}^N\A^k)^c$. Using~\eqref{49}, for all large $n$,
\begin{eqnarray*}
\E\bigl[e^{b_n^2[H^n\wedge\Del_1]}\bigr] &\leq&\sum_{k=1}^N
\E\bigl[e^{b_n^2[H^n\wedge\Del_1]}\chi_{\{(\tilde
{A}^n,\tilde{S}^n)\in\A^k\}}\bigr]+ \E\bigl[e^{b_n^2[H^n\wedge\Del_1]}
\chi_{\{(\tilde{A}^n,\tilde{S}^n)\in
\B\}}\bigr]
\\
&\leq&\sum_{k=1}^N\E
\bigl[e^{b_n^2[\int_0^Th(\varphi^k(s))\,ds
+g(\varphi^k(T))
+\omega(\eps)]}\chi_{\{(\tilde{A}^n,\tilde{S}^n)\in\A^k\}}\bigr]\\
&&{}+ \E\bigl[e^{b_n^2\Del_1}
\chi_{\{(\tilde{A}^n,\tilde{S}^n)\in\B\}}\bigr].
\end{eqnarray*}
Now by Assumption \ref{moderate}, for all large $n$,
\begin{eqnarray*}
\frac{1}{b_n^2}\log\p\bigl(\bigl(\tilde{A}^n,\tilde{S}^n
\bigr) \in\oo{\A^k}\bigr)&\leq&-\inf_{\psi\in\oo{\A^k}}\Ir(\psi)+
\frac
{\eps}{2},\\
 \frac{1}{b_n^2}\log\p\bigl(\bigl(\tilde{A}^n,
\tilde{S}^n\bigr)\in\B\bigr)&\leq&-\inf_{\psi\in\B}\Ir(
\psi)+\eps.
\end{eqnarray*}
Hence for large $n$,
\begin{eqnarray*}
\E\bigl[e^{b_n^2[H^n\wedge\Del_1]}\bigr]&\leq&\sum_{k=1}^N
e^{b_n^2[\int
_0^Th(\varphi^k(s))\,ds
+g(\varphi^k(T))+\omega(\eps)
-\inf_{\psi\in\bar{A}_{v^k}}\Ir(\psi)+{\eps}/{2}]} \\
&&{}+ e^{b_n^2[\Del_1-\inf_{\psi\in\B}\Ir(\psi)+\eps]}
\\
&\leq&\sum_{k=1}^N e^{b_n^2[\int_0^Th(\varphi^k(s))\,ds+g(\varphi
^k(T))-\Ir(\psi^k)
+\omega(\eps)+\eps]} +
e^{b_n^2[\Del_1-\Del+\eps]},
\end{eqnarray*}
where for the first term on the RHS we used (\ref{421}) and for the second
term we used the fact $\B\subset\calD_\Del^c$ and the definition of
$\calD_\Del$.

The last term on \eqref{432} is bounded by
\(
\E[e^{b_n^2(c_{11}\La^n+c_{11}+\La^n-G(\Del_1))}].
\)
From Assumption \ref{unbounded}, there exists a constant $c_{12}$ such
that for all large $n$,
\[
\frac{1}{b_n^2}\log\E\bigl[e^{b_n^2(c_{11}+1)\La^n}\bigr]<c_{12}.
\]
Therefore from (\ref{432}) we obtain
\begin{eqnarray*}
&&\limsup\frac{1}{b_n^2}\log\E\bigl[e^{b_n^2H^n}\bigr]
\\
&&\qquad \leq\max_{1\leq k\leq N} \biggl[\int_0^Th
\bigl(\varphi^k(s)\bigr)\,ds+g\bigl(\varphi ^k(T)\bigr)-\Ir
\bigl(\psi^k\bigr)+\omega(\eps) +\eps \biggr]
\\
&&\qquad\quad{} \vee[\Del_1-\Del+\eps]\vee\bigl[c_{11}+c_{12}-G(
\Del_1)\bigr]
\\
&&\qquad\leq\sup_{\psi\in\calP}\bigl[c\bigl(\psi,\bzeta[\psi]\bigr)+\omega(
\eps )+\eps\bigr]\vee[\Del_1-\Del+\eps] \vee\bigl[c_{11}+c_{12}-G(
\Del_1)\bigr].
\end{eqnarray*}
Now let $\eps\to0$ first, then $\Del\to\infty$, recalling that
$c_{11}$, $c_{12}$ and $G$ do not depend on~$\Del$. Finally let $\Del_1\to\infty$, so $G(\Del_1)\to\infty$,
to obtain
\[
\limsup V^n\bigl(\tilde X^n(0)\bigr)\le \limsup
\frac{1}{b_n^2}\log\E\bigl[e^{b_n^2H^n}\bigr]\leq \sup_{\psi\in\calP}
c\bigl(\psi,\bzeta[\psi]\bigr)=V(x),
\]
where for the first inequality we used \eqref{45}, and for the
equality we used \eqref{305}.
This completes the proof.
\end{pf*}

\begin{remark}[(Relaxed version of Assumption \ref{mini})]
\label{rem3}
We return to Remark~\ref{rem2}(b) regarding a relaxed version of
Assumption \ref{mini},
where continuous minimizers $f_h$ and $f_g$ exist.
Under the relaxed assumption the proof of the lower bound is very
similar to the one we have
presented. As far as the upper bound is concerned,
one can define a policy as in the proof of Theorem \ref{th42}, but with
a jump close to the end of the interval, to account for the fact that
in the solution
of the game, the policy has a jump at $T$ from a point determined by
the minimizer $f_h$ to one determined by $f_g$.
The continuity of the paths $\varphi^k$ is used in the proof of
Theorem~\ref{th42},
and so the modified proof will have to address the jump at the end of
the time interval.
This can be done in a manner similar to the way we treat the jump at
time zero.
However, we do not work out the details here.
\end{remark}

\section{The linear case}\label{sec5}

Section~\ref{sec42} describes a policy for the queueing control problem
that is asymptotically optimal.
While the construction of this policy and its analysis facilitate the
proof of the main result,
they fail to provide a simple, closed-form asymptotically optimal policy.
In this section we focus on cost with either $h$ linear and $g=0$
or $g$ linear and $h=0$, aiming at a simple control policy.
More precisely, the assumption on the functions $h$ and $g$ is slightly
weaker, namely that
%
\begin{equation}
\label{77} h(x)=\sum_{i=1}^Ic_ix_i,\qquad
g(x)=\sum_{i=1}^Id_ix_i,
\end{equation}
where $c_i$ and $d_i$ are nonnegative constants, and, in addition,
%
\begin{equation}
\label{78} c^1\mu^1\geq c^2
\mu^2\geq\cdots\geq c^I\mu^I \quad\mbox{and}\quad
d^1\mu^1\geq d^2\mu^2\geq\cdots
\geq d^I\mu^I.
\end{equation}

We consider the so-called $c\mu$ rule, namely the policy that
prioritizes according to
the ordering of the class labels, with highest priority to class 1. Let
us construct
this policy rigorously by considering the set of equations
%
\begin{eqnarray}
\label{52} B^n_1(t)&=&\chi_{\{X^n_i(t)>0\}},
\nonumber
\\[-8pt]
\\[-8pt]
\nonumber
 B^n_2(t)&=&
\chi_{\{X^n_1(t)=0,X^n_2(t)>0\}},\ldots,B^n_I(t) =
\chi_{\{X^n_1(t)=0,\ldots,X^n_{I-1}(t)=0,X^n_I(t)>0\}}.
\end{eqnarray}
Arguing as in Section~\ref{sec42}, considering \eqref{52}
along with the model equations \eqref{01}--\eqref{03}, it is easy to
see that
there exists a unique solution, this solution is used to define the processes
$X^n,D^n,T^n,B^n$, and moreover $B^n$ is an admissible policy.

The result below states that the policy is asymptotically optimal.
%
\begin{theorem}\label{th5}
Let Assumptions \ref{moderate} and \ref{unbounded} hold and assume
$g$ and $h$
satisfy~\eqref{77}--\eqref{78}.
Then, under the priority policy $\{B^n\}$ of \eqref{52},
\[
\lim_{n\to\infty}J^n\bigl(\tilde{X}^n(0),B^n
\bigr)=V(x).
\]
\end{theorem}

\begin{pf}
As explained in Example \ref{mini-rem},
Assumption \ref{mini} holds.
As a result, the lower bound stated in Theorem \ref{th41} is valid.
It therefore suffices to prove that $\limsup_{n\to\infty}J^n(\tilde
{X}^n(0),B^n)\le V(x)$.
The general strategy of the proof of Theorem \ref{th42} is repeated here;
the details of proving the main estimates are, of course, different.

Thus, given constants $\Del$ and $\eps$ we consider $\calD_\Del$ of
\eqref{60}, $M$, the constants $\delta_1,\delta,\eta,v_0, v_2$, the members
$\psi^k$ of $\calD_\Del$, the sets $\calA^k=\calA_{v^k}(\psi^k)$
and the events $\Om^n_k$ [see \eqref{2222}]
precisely as in the proof of Theorem \ref{th42}.
We also set
\(
(\varphi^k,\xi^k,\zeta^k, w^k)=(\bph[\psi^k],\bxi[\psi^k],\bzeta
[\psi^k],\bw[\psi^k])
\)
as in that proof.

In what follows, $c_1$, $c_2,\ldots$ denote constants independent of
$\Del$, $\eps$, $\delta_1$, $\delta$,
$\eta$, $v_0, v_2$ and $n$.
Analogously to \eqref{41} and \eqref{42}, we aim at proving that
there exists
a constant~$c_1$, such that for all sufficiently large $n$,
%
\begin{equation}
\label{54} \bigl\|\tilde X^n\bigr\|^*_T\le c_1
\bigl(1+\La^n\bigr),
\end{equation}
(where, as before, $\La^n=\|\tilde A^n\|^*_T+\|\tilde S^n\|^*_T$), and
%
\begin{equation}
\label{55} \sup_{[v_2,T]}\bigl\|\tilde X^n-
\varphi^k\bigr\| \le c_1\eps\qquad \mbox{on } \Om
^n_k, k=1,2,\ldots,N.
\end{equation}
Once these estimates are established, the proof can be completed
exactly as in step~5 of
the proof of Theorem \ref{th42}. We therefore turn to proving \eqref
{54} and~\eqref{55}.

Recall that $\theta^n=(\frac{n}{\mu^n_1},\frac{n}{\mu^n_2},\ldots
,\frac{n}{\mu^n_I})$.
Moreover, by \eqref{52}, $\sum B^n_i=0$ holds if and only if for all
$i$, $X^n_i=0$,
equivalently $\theta^n\cdot\tilde X^n=0$.
Therefore by \eqref{53},
\[
\theta^n\cdot Z^n(t)=\frac{\sqrt{n}}{b_n} \Biggl(t-\int
_0^t\sum_{i=1}^IB^n_i(s)\,ds
\Biggr)= \frac{\sqrt{n}}{b_n}\int_0^t
\chi_{\{\theta^n\cdot\tilde X^n(s)=0\}}\,ds.
\]
Hence from (\ref{eqn2}), with
%
\begin{equation}
\label{57} Y^n_{\#,i}(t)=\tilde{X}^n_i(0)
+ y^n_it +\tilde{A}^n_i(t)-
\tilde{S}^n_i\bigl(T^n_i(t)
\bigr),
\end{equation}
we have
%
\begin{equation}
\label{555} \theta^n\cdot\tilde{X}^n(t)=
\theta^n\cdot Y^n_\# + \frac{\sqrt{n}}{b_n}\int
_0^t\chi_{\{\theta^n\cdot\tilde
X^n(s)=0\}}\,ds.
\end{equation}
Since $\theta^n\cdot\tilde X^n$ is nonnegative and $\theta^n\cdot
Z^n$ increases only when
$\theta^n\cdot\tilde X^n$ vanishes, it follows that $(\theta^n\cdot
\tilde X^n,\theta^n\cdot Z^n)$
solve the Skorohod problem for $\theta^n\cdot Y^n_\#$.
As a result,
\[
\bigl|\theta^n\cdot\tilde{X}^n\bigr|^*_T+\bigl|
\theta^n\cdot Z^n\bigr|^*_T \le4 \bigl|
\theta^n\cdot Y^n_\#\bigr|^*_T.
\]
Also, using \eqref{eqn2}, the nonnegativity of $\tilde X^n_i$ implies
\[
Z^n_i(t)\geq- Y^n_{\#,i}(t).
\]
Since $\theta^n\to\theta, y^n_i\to y_i, \tilde X^n(0)\to x$, it
follows that
there exists a constant $c_1$ such that for all $n$ large, \eqref{54}
holds, as well as
%
\begin{equation}
\bigl\|Z^n\bigr\|^*_T \leq c_1\bigl(1+
\La^n\bigr). \label{554}
\end{equation}

Toward proving \eqref{55}, let us compute the paths $\varphi^k$.
As mentioned in Example \ref{mini-rem}, the corresponding minimizing curve
is given by $f(w)=(0,\ldots,0, w\mu^I)$, $w\ge0$. Recall notation
\eqref{70}
and that $\xi^k=\bxi[\psi^k]$. Thus
%
\begin{equation}
\label{58} \varphi^k_i=\cases{ 0, &\quad $\mbox{if }
i=1,2,\ldots,I-1$,\vspace*{2pt}
\cr
\mu^I\Gam\bigl[\theta\cdot
\xi^k\bigr], &\quad $\mbox{if } i=I$.}
\end{equation}

Define $\I^\prime=\{1,2,\ldots, I-1\}$ and $\rho'=\sum_{i=1}^{I-1}\rho_i$.
Then by \eqref{eqn2} and \eqref{06},
\begin{eqnarray*}
\tilde X^{n,\prime}(t)&:=&\sum_{i\in\I^\prime}
\theta^n_i\tilde X^n_i(t)=
\sum_{i\in\I^\prime}\theta^n_i
Y^n_{\#,i}(t) +\frac{\sqrt{n}}{b_n}\sum
_{i\in\I^\prime}\bigl(\rho_it-T^n_i(t)
\bigr)
\\
&=& U^n(t)+\frac{\sqrt{n}}{b_n}\int_0^t
\chi_{\{\tilde X^{n,\prime
}(s)=0\}}\,ds,
\end{eqnarray*}
where
\[
U^n(t)=\sum_{i\in\I^\prime}\theta^n_i
Y^n_{\#,i}(t)+\frac{\sqrt {n}}{b_n}\bigl(\rho'-1
\bigr)t,
\]
and we used \eqref{52} by which $\sum_{\calI'}B^n_i=0$ holds if and
only if
$X^n_i=0$ for all $i\in\calI'$.
Hence, invoking again the Skorohod problem,
%
\begin{equation}
\label{507} \tilde X^{n,\prime}(t)= U^n(t)+\sup
_{s\in[0, t]}\bigl\{-U^n(s)\vee0\bigr\}.
\end{equation}

We will argue that, on $\Om^n:=\bigcup_k\Om^n_k$, for all sufficiently
large $n$,
%
\begin{equation}
\label{56} \sup_{[v_2,T]}\bigl|\tilde X'^n\bigr|
\le c_2\eps.
\end{equation}
To this end, let us fisrt show that, for all sufficiently large $n$,
the following holds:
On $\Om^n$, $U^n(t_2)\leq U^n(t_1)$ whenever\vadjust{\goodbreak} $t_1,t_2\in[0,T]$ are
such that
$t_2-t_1\geq v_2$. Suppose this claim is false. Then there are
infinitely many $n$
for which there exist ($n$-dependent) $t_1,t_2\in[0,T]$ with
$t_2-t_1\geq v_2$ but
$U^n(t_2)> U^n(t_1)$ on $\Om^n$.
Thus
\begin{eqnarray*}
&& \sum_{i\in\I^\prime}\theta^n_i
\bigl[ \tilde{X}^n_i(0) + y^n_it_1+
\tilde{A}^n_i(t_1)-\tilde{S}^n_i
\bigl(T^n_i(t_1)\bigr)\bigr]
\\
&&\quad{}-\sum_{i\in\I^\prime}\theta^n_i
\bigl[ \tilde{X}^n_i(0) + y^n_it_2+
\tilde{A}^n_i(t_2)-\tilde{S}^n_i
\bigl(T^n_i(t_2)\bigr)\bigr]
\\
&&\qquad < \frac{\sqrt{n}}{b_n}\bigl(\rho'-1\bigr) (t_2-t_1)
\leq \frac{\sqrt{n}}{b_n}\bigl(\rho'-1\bigr)v_2.
\end{eqnarray*}
However, this is a contradiction because
the RHS tends to $-\iy$ as $n\to\iy$ whereas the LHS remains bounded.
This proves the claim.

Next, note that, for a similar reason, for all sufficiently large $n$,
$U^n(t)<0$ on~$\Om^n$,
for $t\geq v_2$. Hence for $t\geq v_2$ and $n$ large, we have on $\Om^n$,
\[
\sup_{s\in[0, t]}\bigl\{-U^n(s)\vee0\bigr\}=\sup
_{s\in[0, t]}\bigl\{-U^n(s)\bigr\} =\sup
_{s\in[t-v_2, t]}\bigl\{-U^n(s)\bigr\}.
\]
Thus using (\ref{507}), on $\Om^n$, we have for all $n$ large and
$t\geq v_2$,
%
\begin{eqnarray}
\label{508} \tilde X^{n,\prime}(t) &=& U^n(t)+\sup
_{s\in[t-v_2, t]}\bigl\{-U^n(s)\bigr\}
\nonumber
\\
&\leq& \sum_{i\in\I^\prime}\theta^n_iY^n_{\#,i}(t)
+\frac{\sqrt{n}}{b_n}\bigl(\rho'-1\bigr)t \nonumber\\
&&{}+ \sup_{[t-v_2, t]}
\biggl[-\sum_{i\in\I^\prime}\theta^n_iY^n_{\#,i}(s)
-\frac{\sqrt{n}}{b_n}\bigl(\rho'-1\bigr)s \biggr]
\\
&\leq& \sum_{i\in\I^\prime}\theta^n_iY^n_{\#,i}(t)
+ \sup_{[t-v_2, t]} \biggl[-\sum_{i\in\I^\prime}
\theta^n_iY^n_{\#
,i}(s) \biggr]
\nonumber
\\
&\leq& c_3\eps+ c_3\bigl[\osc_{v_2} \bigl(
\tilde A^n\bigr)+\osc_{v_2}\bigl(\tilde{S}^n
\bigr)\bigr],\nonumber
\end{eqnarray}
where we used \eqref{57} and the fact that $T^n_i$ are Lipschitz with
constant 1. On~$\Om^n_k$,
%
\begin{equation}
\label{513} \osc_{v_2}\bigl(\tilde A^n\bigr)\leq2\bigl\|
\tilde A^n-\psi^{k,1}\bigr\|^*+\osc _{v_2}\bigl(
\psi^{k,1}\bigr) \le3\eps,
\end{equation}
where we used \eqref{457} and \eqref{456}. Similarly, $\osc
_{v_2}(\tilde S^n)\le3\eps$.
Using this in \eqref{508} gives~\eqref{56}.

Next, recall that $\theta^n\cdot\tilde X^n=\Gam[\theta^n\cdot Y^n_\#]$.
Note by \eqref{58} that $\theta\cdot\varphi^k=\Gam[\theta\cdot
\xi^k]$.
Therefore using the Lipschitz property of $\Gam$ we have, for all sufficiently
large~$n$,
%
\begin{eqnarray}
\label{511}
\nonumber
\bigl|\theta^n\cdot\tilde X^n-\theta
\cdot\varphi^k\bigr|^*_T &\leq& 2\bigl|\theta^n\cdot
Y^n_\#-\theta^n\cdot\xi^k\bigr|^*_T +
2\bigl\|\theta^n-\theta\bigr\|\bigl\|\xi^k\bigr\|^*_T
\\
&\le& c_4\bigl\|Y^n_\#-\xi^k
\bigr\|^*_T+\eps
\\
&\le& c_4\sum_i\bigl\{\bigl|\tilde
A^n_i-\psi^{k,1}_i\bigr|^*_T
+\bigl|\tilde S^n_i\circ T^n_i-\brho
\bigl[\psi^{k,2}_i\bigr]\bigr|^*_T\bigr\}+2\eps.\nonumber
\end{eqnarray}
Now, on $\Om^n_k$, $\|\tilde A^n-\psi^{k,1}\|\le\eps$
and $\|\tilde S^n-\psi^{k,2}\|\le\eps$. Moreover, from (\ref{554}),
\[
\sup_{[0,T]}\bigl|\bigl(\rho_i t-T^n_i(t)
\bigr)\bigr|\leq v_2,
\]
on $\Om^n$.
It follows that, on $\Om^n_k$, for all sufficiently large $n$,
%
\begin{equation}
\label{512} \bigl|\theta^n\cdot\tilde X^n-\theta\cdot
\varphi^k\bigr|^*_T \leq c_5\eps+
\osc_{v_0}\bigl(\psi^{k,2}\bigr)\le c_6\eps,
\end{equation}
where the last inequality follows from \eqref{456}.

Now, by \eqref{56} and the fact that $\varphi^k_i=0$ for $i<I$ [see
\eqref{58}], we have
$\sup_{[v_2,T]}|\tilde X^n_i-\varphi^k_i|\le c_7\eps$ for $i<I$, on
$\Om^n_k$ for large $n$.
Combining this with \eqref{512}, the convergence $\theta^n\to\theta
$ and the fact
that the $I$ vectors $\theta$ and $\{e_i, i<I\}$ are linearly
independent, gives
$\sup_{[v_2,T]}\|\tilde X^n-\varphi^k\|\le c_8\eps$, on $\Om^n_k$,
for all sufficiently large $n$.
This proves \eqref{55} and completes the proof of the result.
\end{pf}

\begin{appendix}\label{app}
\section*{Appendix}
\begin{pf*}{Proof of Proposition \ref{prop-bd}}We borrow some
ideas from the proof of
Lemma A.1 in \cite{Puhal-1999}.
Clearly, the statements regarding $\tilde A^n$ and $\tilde S^n$ are
identical, hence it suffices to consider
only the former.
Define $M^i_A(u)=\E[e^{u\IA_i}]$ for $u\in\R$.
It suffices to prove that for any positive $K>0$ and $i\in\I$,
\begin{eqnarray*}
\limsup\frac{1}{b_n^2}\log\E\bigl[e^{b_n^2(K|\tilde{A}^n_i|^*)}\bigr] &< & \infty.
\end{eqnarray*}
Assume $i=1$. Since $M^1_A(u)=\E[e^{u\IA_1}]$
is finite around $0$, it is $C^2$ there, and so is $H^1_A(u):=\log M^1_A(u)$.
Therefore by Taylor expansion there exist $\gamma,\delta>0$ such that
%
\begin{equation}
\label{ap-1} \bigl|H^1_A(u)-u\bigr|\leq\gamma u^2\qquad
\mbox{for all } u \mbox{ with } |u|\leq\delta.
\end{equation}
Here we have used the fact that $\frac{dM^1_A}{du}(0)=\E[\IA_1]=1$.
Note that
\begin{eqnarray*}
&&\E\bigl[ e^{b_n^2(K|\tilde{A}^n_1|^*)}\bigr]
\\
&&\qquad=1+b_n^2K\int_0^\infty
e^{b_n^2Kt}\p\bigl(\bigl|\tilde{A}^n_1\bigr|^*>t\bigr)\,dt\leq
1+b_n^2Ke^{Kb_n^2}\\
&&\qquad\quad{}+b_n^2K
\int_1^\infty e^{b_n^2Kt}\p\bigl(\bigl|
\tilde{A}^n_1\bigr|^*>t\bigr)\,dt.
\end{eqnarray*}
For $t\geq1$,
\begin{eqnarray*}
\p\bigl(\bigl|\tilde{A}^n_i\bigr|^*>t\bigr)&=& \p\bigl(\exists v
\in[0,T] \mbox{ such that } \bigl|\tilde{A}^n_1(v)\bigr|>t\bigr)
\\
&\leq& \p\bigl(\exists v\in[0,T] \mbox{ such that } \tilde {A}^n_1(v)<-t
\bigr)\\
&&{} + \p\bigl(\exists v\in[0,T] \mbox{ such that } \tilde {A}^n_1(v)>t
\bigr).
\end{eqnarray*}
Now
\begin{eqnarray*}
\tilde{A}^n_1(v)>t \quad& \Leftrightarrow&\quad
A^n_1(v)>b_n\sqrt{n}t + \lambda^n_1v,
\\
\tilde{A}^n_1(v)<-t \quad& \Leftrightarrow&\quad
A^n_1(v)<-b_n\sqrt{n}t +
\lambda^n_1v.
\end{eqnarray*}
Let $\lfloor x\rfloor$ denote the largest integer less than
or equal to $x$. Also assume $-b_n\sqrt{n}t + \lambda^n_1T> 0$. Then
\begin{eqnarray*}
&& \p\bigl(\exists v\in[0,T] \mbox{ such that } \tilde{A}^n_1(v)<-t
\bigr)
\\
&&\qquad= \p\bigl(\exists v\in[0,T] \mbox{ such that } A^n_1(v)<-b_n
\sqrt{n}t + \lambda^n_1v\bigr)
\\
&&\qquad\leq \p\Biggl(\exists v\in[0,T] \mbox{ such that } \sum
^{\lfloor
-b_n\sqrt{n}t + \lambda^n_1v+1\rfloor}_{l=1}\IA_1(l)>\lambda^n_1
v\Biggr)
\\
&&\qquad\leq \p\Biggl(\exists v\in[0,T] \mbox{ such that}\\
&& \hspace*{46pt}\sum
^{\lfloor
-b_n\sqrt{n}t + \lambda^n_1v+1\rfloor}_{l=1}\bigl(\IA_1(l)-1\bigr)>\lambda
^n_1 v-\bigl\lfloor-b_n\sqrt{n}t +
\lambda^n_1v+1\bigr\rfloor\Biggr)
\\
&&\qquad\leq \p\Biggl(\exists v\in[0,T] \mbox{ such that } \sum
^{\lfloor
-b_n\sqrt{n}t + \lambda^n_1v+1\rfloor}_{l=1}\bigl(\IA_1(l)-1\bigr)>
b_n\sqrt{n}t-1\Biggr).
\end{eqnarray*}
We define $V_k=\sum_{l=1}^k(\IA_1(l)-1)$. Then $\{V_k\}$ is a
martingale w.r.t. the filtration generated by $\{\IA_1(l)\}$. For all
large $n$ and $t\ge1$, $b_n\sqrt{n}t-1>0$.
Denote $L^n=\lfloor-b_n\sqrt{n}t + \lambda^n_1T+1\rfloor$. Then
\begin{eqnarray*}
\p\bigl(\exists v\in[0,T] \mbox{ such that } \tilde{A}^n_1(v)<-t
\bigr) &\leq& \p\Bigl(\sup_{1\leq k\leq L^n}|V_k|>b_n
\sqrt{n}t-1\Bigr)
\\
&\leq& e^{-\beta_n(b_n\sqrt{n}t-1)}\E\Bigl[\sup_{1\leq k\leq
L^n}e^{\beta_n|V_k|}
\Bigr],
\end{eqnarray*}
where $\beta_n>0$ are any constants. We note that $\{e^{\beta
_n|V_k|}\}_k$ is a sub-martingale. Hence by Doob's martingale inequality
\[
\E\Bigl[\sup_{1\leq k\leq L^n}e^{\beta_n|V_k|}\Bigr]\leq\E\Bigl[\sup
_{1\leq
k\leq L^n}e^{2\beta_n|V_k|}\Bigr]^{{1}/{2}} \leq2\E
\bigl[e^{2\beta_n|V_{L^n}|}\bigr]^{{1}/{2}}.
\]
Thus
\begin{eqnarray*}
&&\p\bigl(\exists v\in[0,T] \mbox{ such that } \tilde{A}^n_1(v)<-t
\bigr)
\\
&&\qquad\leq2e^{-\beta_n(b_n\sqrt{n}t-1)}\E\bigl[e^{2\beta_n|V_{L^n}|}\bigr]
^{{1}/{2}}
\\
&&\qquad\leq2e^{-\beta_n(b_n\sqrt{n}t-1)}\bigl[\E\bigl[e^{2\beta_nV_{L^n}}\bigr]+\E
\bigl[e^{-2\beta_nV_{L^n}}\bigr]\bigr]^{{1}/{2}}
\\
&&\qquad\leq2e^{-\beta_n(b_n\sqrt{n}t-1)}\bigl[e^{L^n(H^1_A(2\beta_n)-2\beta_n)} +e^{L^n(H^1_A(-2\beta_n)+2\beta_n)}
\bigr]^{{1}/{2}}.
\end{eqnarray*}
If $2\beta_n\leq\delta$ and $n$ is large enough so that $\frac
{b_n\sqrt{n}t}{2}-1>0$ holds, then using (\ref{ap-1}) we have
\begin{eqnarray*}
\p\bigl(\exists v\in[0,T] \mbox{ such that } \tilde{A}^n_1(v)<-t
\bigr) &\leq & 2\sqrt{2}e^{-\beta_n({b_n\sqrt{n}t}/{2})}e^{4L^n\gamma\beta_n^2}
\\
&\leq& 2\sqrt{2}e^{-\beta_n({b_n\sqrt{n}t}/{2})}e^{4(-b_n\sqrt
{n}t + \lambda^n_1T+1)\gamma\beta_n^2}.
\end{eqnarray*}
Now we choose $\beta_n=\frac{b_n}{\sqrt{n}}(2K+2)$, and we choose
$n_1$ such that for $n\geq n_1$, $2\beta_n\leq\delta$. Then
%
\begin{eqnarray}
\label{ap-2} &&\p\bigl(\exists v\in[0,T] \mbox{ such that } \tilde{A}^n_1(v)<-t
\bigr)
\nonumber
\\[-8pt]
\\[-8pt]
\nonumber
&&\qquad\leq 2\sqrt{2}e^{b_n^216({\lambda^n_1T+1}/{n})\gamma(K+1)^2}e^{-b_n^2(K+1)t}.
\end{eqnarray}
In a similar way we obtain $n_2$ such that for all $n\geq n_2$,
%
\begin{eqnarray}
\label{ap-3} &&\p\bigl(\exists v\in[0,T] \mbox{ such that }
\tilde{A}^n_1(v)>t
\bigr)
\nonumber
\\[-8pt]
\\[-8pt]
\nonumber
&&\qquad\leq 2\sqrt{2}e^{b_n^216({\lambda^n_1T}/{n})\gamma(K+2)^2}e^{-b_n^2(K+1)t}.
\end{eqnarray}
Thus from (\ref{ap-2}) and (\ref{ap-3}) we have constants $n_3,
\gamma_1,\gamma_2$ such that for all $n\geq n_3$,
$\p(|\tilde{A}^n_i|^*>t)\leq\gamma_1e^{b_n^2\gamma_2}e^{-b_n^2(K+1)t}$.
Hence for $n\geq n_3$,
\[
\int_1^\infty e^{b_n^2Kt}\p\bigl(\bigl|
\tilde{A}^n_1\bigr|^*>t\bigr)\,dt\leq\gamma _1e^{b_n^2\gamma_2}
\int_1^\infty e^{-b_n^2t}\,dt=\frac{1}{b_n^2}
\gamma _1e^{b_n^2(\gamma_2-1)}
\]
and
$
\E[e^{b_n^2(K|\tilde{A}^n_1|^*)}]\leq1+b_n^2Ke^{Kb_n^2}+K\gamma
_1e^{b_n^2(\gamma_2-1)}$,
which gives the required estimate.
\end{pf*}
\end{appendix}

\section*{Acknowledgments} We are grateful to the referees for
constructive comments and suggestions
that have much improved the exposition.

%

%



\printaddresses

\end{document}